\documentclass[a4paper,10pt,twoside]{article}
\usepackage[active]{srcltx}

\newcommand{\be}{\begin{enumerate}}
\newcommand{\ee}{\end{enumerate}}
\newcommand{\beqn}{\begin{eqnarray*}}
\newcommand{\eeqn}{\end{eqnarray*}}
\newcommand{\disp}{\displaystyle}

\newcommand{\incl}[1][r]
      {\ar@<-0.2pc>@{^(-}[#1] \ar@<+0.2pc>@{-}[#1]}

\def\N{{\mathbb N}}
\def\P{{\mathbb P}}

\def\Z{{\mathbb Z}}
\def\Ac{{\mathcal A}}
\def\Bc{{\mathcal B}}
\def\Cc{{\mathcal C}}
\def\Fc{{\mathcal F}}

\def\Ic{{\mathcal I}}
\def\Lc{{\mathcal{L}}}
\def\Mc{{\mathcal M}}
\def\Nc{{\mathcal N}}
\def\Oc{{\mathcal O}}
\def\Pc{{\mathcal P}}
\def\Qc{{\mathcal Q}}

\def\Xc{{\mathcal{X}}}
\def\Yc{{\mathcal{Y}}}
\def\Zc{{\mathcal{Z}}}

\def\Ugo{{\mathfrak U}}
\def\gm{{\mathbb G}_{\rm{m}}}

\def\bgm{{\rm B}{\mathbb G}_{\rm{m}}}

\def\fleche{\rightarrow}
\def\flechelongue{\longrightarrow}

\def\ov{\overline}
\newcommand{\cartesien}{\ar@{}[dr]|{\square}}

\def\lind{\lim_{\longrightarrow}}
\def\lpro{\lim_{\longleftarrow}}

\def\Ker{\hbox{\rm Ker}\,}
\def\Coker{\hbox{\rm Coker}\,}

\def\Qcoh{\hbox{${\mathfrak{Qcoh}}$}\,}

\def\Aut{\hbox{\rm Aut}}

\def\Isom{\hbox{\rm Isom}}

\def\Homs{\hbox{\rm \underline{Hom}}}
\def\fIsom{{\mathcal I\!som}\,}

\def\Im{\hbox{\rm Im}\,}

\def\id{\hbox{\rm id}}

\def\Gal{\hbox{\rm Gal}\,}

\def\Hom{\hbox{\rm Hom}}

\def\fHom{{\mathcal{H}om}\,}

\def\red{\text{\rm red}}
\def\Spec{{\rm Spec}\,}
\def\Specc{{\mathcal{S}\!\it{pec}}\,}
\def\Tor{{\rm Tor}\,}

	{\begin{pmatrix}}%
	{\end{pmatrix}}
	{\begin{vmatrix}}%
	{\end{vmatrix}}

	{\ensuremath{\left \{ \hskip -1.5 mm \begin{array}{#1@{\ =\ }l}}}%
	{\end{array}\right.}

	{\ensuremath{\left \{ \hskip -1.5 mm%
	 \begin{array}{#1@{\ }c@{\ }l}}}%
	{\end{array}\right.}

	{\ensuremath{\left \{ \hskip -1.5 mm \begin{array}{#1@{\quad}l}}}%
	{\end{array}\right.}
\newcommand{\fonction}[5]{%
        \ensuremath{#1\colon
        \left\{\hskip -1.5 mm                   
        \begin{array}{c@{\ }c@{\ }l}
        \medskip #2 & \longrightarrow & #3 \\
        #4 & \longmapsto & #5 \\
        \end{array}
        \right .
        }}

\newtheorem{souscor}[subsubsection]{Corollary}

\newtheorem{sousprop}[subsubsection]{Proposition}

\newtheorem{sousthm}[subsubsection]{Theorem}

\newtheorem{sousdefi}[subsubsection]{Definition}

\newtheorem{souslem}[subsubsection]{Lemma}

\newtheorem{sousremarque}[subsubsection]{Remark}

\newtheorem{sousexemple}[subsubsection]{Example}

\newenvironment{demo}{{\bf Proof.\,}}{$\square$ \vskip .3cm}

\newenvironment{etape}[1]{$\bullet$ \emph{#1}\\}{\vskip .2cm}

\usepackage{latexsym}
\usepackage{amsfonts}
\usepackage{amstext}
\usepackage{amsmath}
\usepackage{graphicx}
\usepackage{enumitem}
\setlist{noitemsep,topsep=0pt}
\setitemize{noitemsep,topsep=0pt}
\setenumerate{noitemsep,topsep=0pt}
\usepackage[all,dvips]{xy}
\input xy
\xyoption{all}
\xyoption{web}
\usepackage{amssymb,euscript,multicol,graphicx}
\usepackage[ec]{aeguill}       
\usepackage{delarray, fancyhdr}
\author{Sylvain Brochard\\ I3M, Universit\'e de Montpellier 2\\
\textsf{brochard@math.univ-montp2.fr}}
\SetMathAlphabet\mathcal{normal}{U}{rsfs}{m}{n} 

\renewcommand{\sectionmark}[1]{\markboth{#1}{}}

\fancypagestyle{plain}{ 
\fancyhead{} 
}

\title{Finiteness theorems for the Picard objects of an algebraic stack
}
\date{\empty}

\def\Tor{{\rm Tor}\,}
\def\Pic{\textrm{\rm Pic}}

\def\pic{\Pic_{\X/S}}
\def\pictox{\Pic_{\X/S}^{\tau}}
\def\pictoy{\Pic_{\Y/S}^{\tau}}
\def\picto{\Pic^{\tau}}
\def\piczero{\Pic_{\X/S}^0}
\def\champic{\mathcal{P}ic}

\def\X{\mathcal{X}}
\def\Y{\mathcal{Y}}
\def\L{\mathcal{L}}
\def\Xt{\widetilde{\mathcal{X}}}
\def\nersev{\textrm{\rm NS}}

\renewcommand{\fleche}{%
\xymatrix@C=1pc{\ar[r] &}}
\def\flechelongue{%
\xymatrix{\ar[r] &}}

\entrymodifiers={!!<0pt,0.7ex>+}

\begin{document}

\maketitle


\begin{abstract}
We prove some finiteness theorems for the Picard functor of an algebraic stack, in the spirit of~SGA 6, exp.~XII and~XIII. In particular, we give a stacky version of Raynaud's relative representability theorem, we give sufficient conditions for the existence of the torsion component of the Picard functor, and for the finite generation of the N\'eron-Severi groups or of the Picard group itself. We give some examples and applications. In an appendix, we prove the semicontinuity theorem for a (non-necessarily tame) algebraic stack.

\medskip

\end{abstract}

\tableofcontents

\section{Introduction}

Let $X$ be a proper, flat and finitely presented scheme (or stack) over a base scheme $S$, and let $\Pic_{X/S}$ denote its Picard functor. If $S$ is the spectrum of a field, the connected component of the identity $\Pic_{X/S}^0$ is an open and closed group subscheme. One reason why $\Pic_{X/S}^0$ is a bit easier to handle than the whole Picard functor is that it is of finite type over~$S$. Over a general base scheme $S$, the situation is more complicated and we need additional assumptions to ensure the existence of the neutral component. There are a few positive results: for instance, if $S$ is a characteristic-zero scheme, and if $X \fleche S$ is proper, smooth and with geometrically connected fibers, then $\Pic_{X/S}^0$ is an open subspace of $\Pic_{X/S}$, which is moreover proper over $S$. But without these assumptions, it can actually happen that $\Pic_{X/S}^0$ does not exist as an algebraic space (even if $\Pic_{X/S}$ does), see~\ref{exemple_piczero_non_repres} for an easy example. In such a case, there is an other subfunctor which can serve as a substitute: the torsion component of $\Pic_{X/S}$, denoted by $\picto_{X/S}$. Roughly speaking, it is the set of points of $\Pic_{X/S}$, a power of which lies in $\Pic_{X/S}^0$. In the case of schemes, we know from~SGA~6, exp.~XIII that $\picto_{X/S}$ does exist as an algebraic space much more often than $\Pic^0_{X/S}$, and that it is an open subspace of $\Pic_{X/S}$ which is \emph{of finite type} over $S$. Actually it is in most cases the biggest finite type open group subspace of $\Pic_{X/S}$.

On the other hand, even if $\Pic^0_{X/S}$ exists, a drawback of this subfunctor is that, by construction, it does not contain the discrete information that is enclosed in the whole Picard functor (the number of connected components for instance). For this reason, it is natural to study the N\'eron-Severi groups. For every geometric point $\ov{s}$ of $S$, the N\'eron-Severi group at $\ov{s}$ is
$$\nersev(\ov{s})=\Pic_{X_{\ov{s}}/\kappa(\ov{s})}(\kappa(\ov{s}))/\Pic_{X_{\ov{s}}/\kappa(\ov{s})}^0(\kappa(\ov{s})).$$
Raynaud and Kleiman proved that these groups are of finite type as soon as $X$ is a proper scheme over $S$. Moreover, if $S$ is noetherian, their rank and the order of their torsion subgroup are uniformly bounded over $S$.

The original motivation of the present work was to provide some information about the N\'eron-Severi groups and the torsion component of the Picard functor of an algebraic stack $\X$ over $S$. It turned out that it was necessary to prove a stacky version of the relative representability theorem from SGA~6, exp.~XII. The latter is the main theorem of this article, stated in~\ref{thm_ppal}. The others are easy consequences of this one, using Chow's lemma. By the way, we also get some finiteness results for the Picard stack $\champic(\X/S)$. 

While proving these results, a special care was taken to avoid superfluous tameness assumptions. Because of the fact that an algebraic stack might have infinite cohomological dimension, it is often convenient to assume that the algebraic stacks under consideration are tame. However, in positive characteristic, this assumption is quite restrictive. In the appendix, we will see that the cohomology of an arbitrary algebraic stack is tractable as soon as the base scheme is regular and noetherian. As a byproduct we get the semicontinuity theorem for algebraic stacks~(\ref{semicont}).

\paragraph*{Contents.}
The paper is organized as follows. In section~\ref{proof_main}, we prove the main theorem. We start with some preliminary lemmas that have been separated from the rest of the proof for the sake of clarity. The sections~\ref{cas_immersion_fermee} and~\ref{cas_general} contain the heart of the proof of~\ref{thm_ppal}. At the end of~\ref{cas_general}, we also give a few immediate corollaries. In particular, we prove that the Picard functor of a proper algebraic stack over a field is always a scheme, without any further assumption (\ref{cor_sur_un_corps}).

In section~\ref{Corollaries}, we give some corollaries and applications. We begin with general finiteness results (\S\ref{general}). For instance we prove that the ``$n^{\textrm{th}}$-power'' maps of the Picard functor $\Pic_{X/S}$ and of the Picard stack $\champic(\X/S)$ are of finite type. The paragraph~\ref{prelim} is devoted to the definition and basic properties of quasi-compact or quasi-separated morphisms of non representable functors. (These elementary results were necessary for~\ref{general}.) The results concerning the torsion component $\pictox$ are in paragraph~\ref{torsion}. We conclude this part with some applications that are more arithmetic in nature (\S\ref{arithmetic}), including the finite generation of the N\'eron-Severi groups mentioned above. We also prove there a few finiteness properties for the Picard \emph{group} of an algebraic stack. Some concrete examples and computations can also be found in this section~\ref{Corollaries}.

In the appendix (\S \ref{appendice}), we provide some cohomological stuff, in the spirit of Mumford's~\cite{Mumford_Abelian_Varieties}~\S~5 ``Cohomology and base change''. In particular, we prove there the semicontinuity theorem for algebraic stacks, that was still lacking to the literature.

\paragraph*{Notations and conventions.}
Following~\cite{LMB}, all algebraic stacks (\emph{a fortiori} all schemes and algebraic spaces) are supposed to be quasi-separated. A stack admitting a smooth cover, and the diagonal of which is (only) representable and locally of finite type will be called an \emph{algebraic stack in Artin's sense}. The cohomology groups on an algebraic stack $\X$ are computed with respect to the smooth-\'etale topology. Note that if $\X$ is a Deligne-Mumford stack, and in particular a scheme, we recover the \'etale cohomology groups (\emph{cf.}~\cite{Brochard_Picard}~A.1.6). If $\X$ is an algebraic stack over a base scheme $S$, its Picard functor, denoted by $\pic$, is the \emph{fppf} sheaf associated with the presheaf $U\mapsto \Pic(\X\times_S U)$. Its Picard stack $\champic(\X/S)$ is the stack whose fiber category over an $S$-scheme $U$ is the category of invertible sheaves on $\X\times_S U$.

\section{A relative representability theorem}
\label{proof_main}
The main result of this section is the relative representability theorem~\ref{thm_ppal}, which generalizes~\cite{SGA6} exp.~XII,~1.1. It says the following. Let $f : \X  \fleche \Y$ be a surjective morphism between two proper and finitely presented algebraic stacks over an integral base scheme~$S$. Then, over a nonempty open subset of $S$, the induced morphism $f^*$ between the Picard functors is quasi-affine and of finite type (see~\ref{thm_ppal} for a more complete statement). Our result is more general in two respects: first it also applies to algebraic stacks, second it gives information not only about the Picard functors, but also about the Picard stacks. For the proof, as in~\cite{SGA6}~XII we treat separately the case of a nilpotent immersion, and then we use non-flat descent. However, the proof given here is much simpler than the proof given (for the particular case of schemes) in~\cite{SGA6}, due to the use of the representability theorem by an algebraic space for $\pic$ (see \cite{Global_Analysis_1} and~\cite{Brochard_Picard}). We do not need to prove the representability of $f^*$, but only its quasi-affineness and its quasi-compactness.

\subsection{Preliminary lemmas}
\label{lemmes_prelim}

In the proof of~\ref{thm_ppal} we will need some improvements of the results given in~\cite{Brochard_Picard}. The point~1) of the theorem below is a refinement of~\cite{Brochard_Picard}~1.2. The superfluous assumption that $\xymatrix@C=1pc{\Oc_{S} \ar[r]^-{\sim} &f_*\Oc_{\X}}$ holds universally has been removed. See the notations and conventions above for the definition of the Picard functor and the Picard stack of an algebraic stack.

\begin{sousthm}
\label{champicqsep}
Let $S$ be scheme and let $f : \X \fleche S$ be a proper, flat and finitely presented algebraic stack over $S$.
\begin{itemize}
 \item[(1)] Then the diagonal $\Delta_{\Pc}$ of the Picard stack $\champic(\X/S)$ is separated and of finite presentation. (In particular the stack $\champic(\X/S)$ is quasi-separated, hence is algebraic in the sense of~\cite{LMB}.)
\item[(2)] If $\Oc_{S} \xymatrix@C=1pc{\ar[r]^-{\sim} & f_*\Oc_{\Xc}}$ holds universally, then $\Delta_{\Pc}$ is quasi-affine, and the Picard functor $\pic$ is a (quasi-separated) algebraic space.
\item[(3)] If $S$ is noetherian and $\X$ is tame\footnote{We recall the definition of a tame stack, due to Abramovich, Olsson and Vistoli, in~\ref{champs_moderes}. In particular, a scheme or an algebraic space is tame.}, then $\Delta_{\Pc}$ is affine and of finite presentation.
\item[(4)] If $U$ is a noetherian scheme with finite global cohomological dimension (\emph{e.g.} a regular scheme of finite dimension) then for any morphism $$\varphi : U \flechelongue \champic(\X/S)\times_S \champic(\X/S),$$ the pullback of $\Delta_{\Pc}$ along $\varphi$ is affine and of finite presentation.
\item[(5)] If $S$ is integral and noetherian, then there is a nonempty open subset~$U$ of~$S$ over which~$\Delta_{\Pc}$ is quasi-affine.
\end{itemize}
\end{sousthm}
\begin{sousremarque} \rm
In view of~(3), (4) and~(5), it is tempting to believe that, under the assumptions of the theorem, the diagonal~$\Delta_{\Pc}$ is always affine. I do not know whether this is the case or not.
\end{sousremarque}
\begin{demo}
 (3) and~(4) are reformulations of~\ref{lem_homisom}. Let us prove~(1). By~\cite{Brochard_Picard} we already know that~$\champic(\X/S)$ is an algebraic stack in Artin's sense. In particular, for a scheme~$U$ and two invertible sheaves~$\Lc, \Mc$ on~$\Xc\times_S U$, the sheaf~$\fIsom(\Lc,\Mc)$ is an algebraic space locally of finite presentation. To prove that it is separated, we can assume that $S$ is the spectrum of a discrete valuation ring and then it is a consequence of (4). It remains to prove that it is quasi-compact. By standard limit arguments (use~\cite{Olsson_Hom_stacks}~2.2) we can assume that $S$ is of finite type over $\Spec \Z$. Moreover we can assume that it is reduced. By noetherian induction, it is enough to prove the assertion over a nonempty open subset $U$ of $S$. Hence we can assume that $S$ is regular and once again it is a consequence of (4).

Let us prove~(2). The assertion about $\pic$ was proved in~\cite{Brochard_Picard}~2.3.3. Now notice that the diagonal of $\pic \fleche S$ is a monomorphism of finite type of algebraic spaces, hence, applying~\cite{LMB}~(A.2.2) it is quasi-affine. Thus to prove that $\Delta_{\Pc}$ is quasi-affine, it is enough to prove that the diagonal of $\champic(\X/S) \fleche \pic$ is quasi-affine. By faithfully flat descent we may assume that the morphism $f$ has a section. Then by~\cite{Brochard_Picard}~2.3.4 $\champic(\X/S)$ is isomorphic to $\pic \times_S \bgm$. The assertion follows since the diagonal of~$\bgm$ is affine.

Finally let us prove (5). We can assume that $S$ is affine and regular. Using~\ref{kunneth_tame} and generic flatness, we can assume that $f$ is cohomologically flat in dimension zero (\emph{i.e.} that forming $f_*\Oc_{\Xc}$ commutes with base change). Moreover we can assume that $f_*\Oc_{\Xc}$ is free of finite rank. Let $\ov{S}=\Spec (f_*\Oc_{\Xc})$ and let $$\xymatrix{\Xc \ar[r]^{\ov{f}} & \ov{S} \ar[r]^{h} & S}$$
be the Stein factorization. There is a dense open subset $\ov{U}$ of $\ov{S}$ over which~$\ov{f}$ is flat. Since $h$ is flat, this dense open subset contains the generic fiber of~$h$. Now, since $h$ is finite, $h(\ov{S}\setminus \ov{U})$ is a closed subset of $S$, that does not contain the generic point of $S$. Thus, replacing $S$ with a nonempty open subset, we may assume that $\ov{f}$ is flat. Owing to~(2), we then know that the diagonal of~$\champic(\Xc/\ov{S}) \fleche \ov{S}$ is quasi-affine. Let $T$ be an $S$-scheme and let $\Lc, \Mc$ be invertible sheaves on $\Xc\times_S T$. We have to prove that $I:=\fIsom_{T}(\Lc, \Mc)$ is quasi-affine over $T$. Let $\ov{T}=\ov{S}\times_S T$. By the above we know that~$\ov{I}:=\fIsom_{\ov{T}}(\Lc, \Mc)$ is quasi-affine. But $I$ is the Weil restriction of $\ov{I}$ along the morphism $h$, and it is known that Weil restriction along a finite and locally free morphism preserves quasi-affineness (see \cite{BRL}~\S 7.6). This concludes the proof.
\end{demo}

The following lemma will be strengthened in~\ref{cor_repres}, using~\ref{thm_ppal}.

\begin{souslem}
\label{lem_pic}
 Let $S$ be an integral scheme and $\X$ a proper and finitely presented algebraic stack over $S$. Then there is a nonempty open subset $U$ of $S$ such that $(\pic)_{|_U}$ is an algebraic space locally of finite presentation over $U$.
\end{souslem}
\begin{demo}
 We can assume that $S$ is affine and of finite type over $\Spec \Z$. Replacing $S$ with a nonempty open subset, we can assume that $S$ is regular, hence has finite global dimension. Moreover by generic flatness (\cite{EGA} EGA~$\textrm{IV}_2$~6.9.1), we can assume that $\X$ is flat over $S$. Then by~\ref{lem_homisom}, the functor $I : T \mapsto \Gamma(\X_T, \Oc_{\X_T})^{\times}$ is an affine scheme of finite presentation over $S$. Note that forming $I$ commutes with base change. Using generic flatness again, we can assume that $I$ is flat over $S$. Now, since we know by~\ref{champicqsep} that $\champic(\X/S)$ is an algebraic stack in the sense of~\cite{LMB}, we can apply~\cite{LMB}~(10.8) and this gives the result. 
\end{demo}

The lemma below is a slight variation of~\cite{LMB}~(A.2.2). It will be one of the key ingredients of the proof of~\ref{thm_ppal}.

\begin{souslem}
\label{lem_quaff}
 Let $f : X \fleche Y$ be a monomorphism of algebraic spaces. Assume that $Y$ is locally noetherian and that $f$ is locally of finite type. Then $f$ is quasi-affine.
\end{souslem}
\begin{demo}
 Using \cite{LMB}~(A.2.2), it is enough to prove that $f$ is quasi-compact. We can assume that $Y$ is affine and reduced. Moreover, using noetherian induction on $Y$, it is enough to prove that there is a nonempty open subset $V$ of $Y$ such that $X\times_Y V$ is quasi-compact. Now we can assume that $Y$ is integral, and using generic flatness (\cite{EGA} EGA~$\textrm{IV}_2$~6.9.1), that $f$ is flat on a nonempty quasi-compact open subspace $U$ of $X$. But now $f$ is open on $U$. Let $V$ denote the image of $U$ in $Y$. Then, since $f$ is a monomorphism, the subspace $f^{-1}(V)$ is equal to $U$ and thus quasi-compact. This concludes the proof.
\end{demo}

Let $S$ be a scheme and let $f : \X \fleche \Y$ be a morphism of algebraic stacks over~$S$. Let $g$ denote the canonical morphism from $\X\times_{\Y} \X$ to $\Y$. We will say that $P(f)$ holds if the diagram
\label{pagepdef}
$$\xymatrix{\Oc_{\Y} \ar[r] & f_*\Oc_{\X} \ar@<0.5ex>[r]^-{\pi_1^*} \ar@<-0.5ex>[r]_-{\pi_2^*} &g_*\Oc_{\X\times_{\Y}\X}}$$
is exact. We will say that $P(f)$ holds universally if $P(f\times_SS')$ holds for any $S$-scheme $S'$. The following lemma generalizes~\cite{SGA6}~XII~2.6. It will allow us to use descent techniques.

\begin{souslem}
\label{lem_factorisation}
 Let $f : \X \fleche \Y$ be a finite morphism of algebraic stacks. Assume that $\Y$ is noetherian, and that $\Oc_{\Y} \fleche f_*\Oc_{\X}$ is injective (\emph{e.g.} $f$ is surjective and $\Y$ is reduced). Then there is a factorization of $f$ as follows
$$\xymatrix{\X=\Y_0 \ar[r]^-{f_1} & \Y_1 \ar[r]^{f_2}& \dots \ar[r]^-{f_n} & \Y_n=\Y}$$
where for each $i$, the morphism $f_i$ is finite and $P(f_i)$ holds.
\end{souslem}
\begin{demo}
 Let $Y\fleche \Y$ be a presentation of $\Y$ with $Y$ affine and let us take the following notations
\beqn
 \Ac &=& \Oc_{\Y}\\
\Bc_0 &=& f_*\Oc_{\X}\\
\Bc_i &=& \xymatrix{\Ker(\Bc_{i-1} \ar@<0.5ex>[r] \ar@<-0.5ex>[r] &\Bc_{i-1}\otimes_{\Ac} \Bc_{i-1})}\\
\Y_i &=& \Specc \Bc_i \quad \quad \textrm{(in the sense of \cite{LMB} 14.2)}\\
Y_i &=& \Y_i \times_{\Y} Y\\
X &=& \X \times_{\Y} Y
\eeqn
Since the morphism $\Y_i\fleche \Y$ is affine, $Y_i$ is an affine scheme, say $Y_i=\Spec B_i$. Moreover, since forming a pushforward or a kernel commutes with flat base change, we have $B_0=f_*\Oc_X$ and
$$B_i = \xymatrix{\Ker(B_{i-1} \ar@<0.5ex>[r] \ar@<-0.5ex>[r] &B_{i-1}\otimes_{A} B_{i-1})}$$
where $A=\Oc_Y$. The proof given in \cite{SGA6}~XII~2.6 shows that for all $i$, the property $P(Y_i\fleche Y_{i+1})$ holds and that for $i$ big enough, $Y_i\fleche Y$ is an isomorphism. By faithfully flat descent, the same holds for the $\Y_i$'s. 
\end{demo}

\begin{souslem}
 \label{fgregc}
Let $S$ be an integral scheme and $\X$ a proper and finitely presented algebraic stack over $S$. Then there is a nonempty open subset $V$ of $S$, an integral scheme $V'$, and a finite flat morphism $V'\fleche V$ such that the fibers of $(\X_{V'})_{\red}$ are geometrically reduced, and the connected components of $\X_{V'}$ have geometrically connected fibers.
\end{souslem}
\begin{demo}
 By standard arguments we can assume that $S=\Spec A$ with $A$ of finite type over $\Z$. Replacing $S$ with a nonempty open subset, we can further assume that $\X$ is flat and (using~\ref{kunneth_tame}) that forming $H^0(\X,\Oc_{\X})$ commutes with any base change. Let $X$ be the spectrum of $H^0(\X,\Oc_{\X})$. We observe that for any field $k$ and any morphism $\Spec k \fleche S$, the stack $\X\times_S \Spec k$ is connected if and only if $X\times_S \Spec k$ is. We will use this (and the fact that being reduced is a property of local nature for the smooth topology) to apply some results from EGA~\cite{EGA} to the stack $\X$. Let $\eta$ be the generic point of $S$, and $k$ its function field. By EGA~$\textrm{IV}_2$~4.6.8, there is a finite extension $k'$ of $k$ such that $(\X_{k'})_{\red}$ is geometrically reduced, and such that the connected components of $\X_{k'}$ are geometrically connected. There is a finite, flat and integral $A$-algebra $A'$ with function field equal to $k'$. Replacing $A$ by $A'$ we can assume that $k'=k$. Now the generic fiber of $\X_{\red}$ is geometrically reduced, hence by EGA~$\textrm{IV}_3$~9.7.7 there is a nonempty open subset $U$ of $S$ over which the fibers of~$\X_{\red}$ are geometrically reduced. We replace $S$ with $U$. Let us denote by $\X_i$, $i=1,\dots, n$ the irreducible components of $\X$. Shrinking $S$ if necessary, we can assume that $\X_i\cap \X_j$ is empty if and only if $(\X_i)_{\eta}\cap(\X_j)_{\eta}$ is empty. We can also assume that $\X$ is connected. Now $\X_{\eta}$ is connected as well, hence geometrically connected by the above. We apply EGA~$\textrm{IV}_3$~9.7.7 to find a nonempty open subset of $S$ over which the fibers of $\X$ are geometrically connected, and this concludes the proof.
\end{demo}

\subsection{The case of a nilpotent immersion}
\label{cas_immersion_fermee}

\begin{sousthm}
\label{cas_nilpotent}
 Let $S$ be an integral scheme and let $f : \X \fleche \Y$ be a surjective closed immersion of proper and finitely presented algebraic stacks over $S$. Then there is a nonempty open subset $U$ of $S$ with the following properties:
\begin{itemize}
 \item[a)] The functors $(\pic)_{|_U}$ and $(\Pic_{\Y/S})_{|_U}$ are algebraic spaces, and the morphism $f^* : (\Pic_{\Y/S})_{|_U} \fleche (\pic)_{|_U}$ is affine and of finite type.
\item[b)] The stacks $\champic(\X/S)_{|_U}$ and $\champic(\Y/S)_{|_U}$ are algebraic stacks, and the morphism $f^* : \champic(\Y/S)_{|_U} \fleche \champic(\X/S)_{|_U}$ is of finite type with affine diagonal.
\end{itemize}
\end{sousthm}
\begin{demo}
We can assume that $S$ is affine, of finite type over $\Z$, and regular. The ideal~$\Ic$ in $\Y$ defining $f$ is a nilpotent ideal. Obviously we can assume that $\Ic$ is square-zero. Using~\ref{lem_pic} we can assume that $\pic$ and $\Pic_{\Y/S}$ are algebraic spaces locally of finite presentation. Let~$k$ be the global dimension of the ring~$\Oc_S$. Using generic flatness we can assume that $\X$, $\Y$, $\Ic$ and all the sheaves $R^i\mu_*\Ic$, $i\leq k+2$ are flat over $S$\,, where~$\mu$ is the structural morphism $\Y\fleche S$. Then by~\ref{champicqsep} the Picard stacks are algebraic, and by~\ref{images_directes_superieures_affines} the sheaves $V_i$ on $(\text{\rm Sch}/S)^{\circ}$ defined on affine schemes by
$$V_i(T)=H^i(\X_T, \Ic_T)$$
are affine schemes of finite type over $S$ for $i=0, 1, 2$.

For each $S$-scheme $T$, the exact sequence of abelian sheaves
$$\xymatrix{0\ar[r] &\Ic_T \ar[r] &\Oc_{\Y_T}^{\times} \ar[r] & (f_T)_*\Oc_{\X_T}^{\times} \ar[r] &0}$$
induces the following exact sequence of abelian groups
$$\xymatrix{0\ar[r] &H^1(\Y_T,\Ic_T) \ar[r] &\Pic(\Y_T) \ar[r] & \Pic(\X_T) \ar[r]^-{\omega} &H^2(\Y_T,\Ic_T)} \quad (*)$$
(see~\cite{Brochard_Picard}~3.1.3 for the details of this computation). Applying the functor ``\emph{fppf} associated sheaf'', we get an exact sequence of \emph{fppf} abelian sheaves:
$$\xymatrix{0\ar[r] &V_1 \ar[r] &\Pic_{\Y/S} \ar[r] & \pic \ar[r] &V_2.}$$
Let $P$ denote the kernel of the map $\pic \fleche V_2$. Since $V_2$ is separated, the morphism $P\fleche \pic$ is a closed immersion thus an affine morphism. Let us prove that the morphism $\Pic_{\Y/S}  \fleche P$ is affine. Since the sequence
$$\xymatrix{0\ar[r] &V_1 \ar[r] &\Pic_{\Y/S} \ar[r] & P \ar[r] &0}$$
is exact, we see that $\Pic_{\Y/S}$ is a $V_1$ pseudo-torsor over $P$. Considering a presentation of the algebraic space $P$ by an affine scheme, we can assume that $P$ is an affine scheme and even that $P=S$. Moreover, by \emph{fppf} descent, we can assume that the morphism $\Pic_{\Y/S} \fleche P$ has a section. Then $\Pic_{\Y/S}$ is a trivial torsor, hence isomorphic to $V_1$ and affine.

Now let us prove b). Let $\Pc$ denote the kernel of the natural map from $\champic(\X/S)$ to $V_2$, \emph{i.e.} $\Pc$ is defined by the cartesian square
$$\xymatrix{\Pc \ar[r] \ar[d] \cartesien & S\ar[d]^e \\ \champic(\X/S) \ar[r] & V_2}$$
where $e : S \fleche V_2$ is the neutral section. Then $\Pc$ is a closed substack of $\champic(\X\!/S)$ and it suffices to prove that the natural morphism from $\champic(\Y\! /S)$ to $\Pc$ is of finite type with affine diagonal. It is enough to prove that the same holds for the morphism obtained after a base change \emph{via} $u : U\fleche \Pc$, where $U$ is an affine scheme. To avoid heavier notations, let us assume that $U=S$ (this is harmless). Now consider the cartesian square
$$\xymatrix{\Qc \ar[r] \ar[d] \cartesien & S\ar[d]^u \\ \champic(\Y/S) \ar[r] & \Pc\, .}$$
By definition of $\Pc$, the point $u$ corresponds to an invertible sheaf $\Lc$ on $\X$, the obstruction class $\omega(\Lc)$ of which is trivial (see~(*), note also that $\omega(\Lc)$ is really trivial and not only locally for the \emph{fppf} topology, because in the \emph{fppf} sheafification process above $V_2$ was not affected). Let us describe the stack $\Qc$. If $T$ is an affine scheme over $S$, the fiber category $\Qc(T)$ is the category of couples $(\Mc,\alpha)$ where $\Mc$ is an invertible sheaf on $\Y\times_S T$
and $\alpha$ is an isomorphism $\alpha : f^*\Mc \fleche \Lc_{|_T}$.
Since the obstruction class $\omega(\Lc)$ is zero, we see from the above exact sequence $(*)$ that the stack $\Qc$ has an $S$-point. Thus we can assume that $\Lc$ is trivial. Now we can check easily that the set of isomorphism classes of objects of $\Qc(T)$ is isomorphic to $\Ker(\Pic(\Y_T) \fleche \Pic(\X_T))$, \emph{i.e.} to $H^1(\Y_T,\Ic_T)=V_1(T)$. Moreover, the group of automorphisms of an object is $\Ker(\Aut(\Oc_{\Y}) \fleche \Aut(\Oc_{\X}))$ \emph{i.e.} it is $H^0(\Y_T,\Ic_T)=V_0(T)$. Hence we see that the algebraic stack $\Qc$ is (isomorphic to) the $S$-groupoid associated to $[\xymatrix@C=1pc{V_0\times_S V_1 \ar@<0.5ex>[r] \ar@<-0.5ex>[r] & V_1}]$ (see~\cite{LMB}~2.4.3), where both maps (source and target) from $V_0\times_S V_1$ to $V_1$ are equal to the projection map on $V_1$. It follows immediately that $\Qc$ is quasi-compact and has an affine diagonal.
\end{demo}

\subsection{General case (non-flat descent arguments)}
\label{cas_general}

Let us start with a lemma that will be helpful when we compare the Picard stack and the Picard functor.

\begin{souslem}
\label{comparaison}
 Let $f : \Xc \fleche \Yc$ be a morphism of algebraic $S$-stacks. Let us denote by $a$ and $b$ the structural morphisms of $\Xc$ and $\Yc$. Assume that $\xymatrix@C=1pc{\Oc_S \ar[r]^-{\sim} & a_*\Oc_{\Xc}}$ and $\xymatrix@C=1pc{\Oc_S \ar[r]^-{\sim} & b_*\Oc_{\Yc}}$ hold universally. Then the natural square of stacks
$$\xymatrix{
\champic(\Yc/S) \ar[r]^{f^*} \ar[d] & \champic(\Xc/S)\ar[d] \\
\Pic_{\Yc/S} \ar[r]^{f^*}& \Pic_{\Xc/S}
}$$
is 2-cartesian.
\end{souslem}
\begin{demo}
\def\FP{\textrm{FP}}
 Let us denote by $\FP$ the fiber product $\Pic_{\Yc/S}\times_{\Pic_{\Xc/S}} \champic(\Xc/S)$. For any $S$-scheme $U$, the category $\FP(U)$ can be described as follows. Its objects are couples $(m,\Lc)$ where $m\in \Pic_{\Yc/S}(U)$ and $\Lc$ is an invertible sheaf on $\Xc\times_S U$, such that the class $[\Lc]$ of $\Lc$ in $\Pic_{\Xc/S}(U)$ is equal to the inverse image $f^*m$. If $(m_1, \Lc_1)$ and $(m_2,\Lc_2)$ are two such pairs, then the set $\Hom((m_1, \Lc_1),(m_2,\Lc_2))$ of morphisms in $\FP(U)$ is empty if $m_1\neq m_2$, and is $\Isom(\Lc_1, \Lc_2)$ if $m_1=m_2$. We have to prove that the natural map $\psi$ from $\champic(\Yc/S)$ to $\FP$ induced by the above square is an isomorphism. Owing to \cite{LMB}~(3.7.1), it suffices to prove that it is a monomorphism and an epimorphism.

Let $(m,\Lc)\in \FP(U)$. Then, \'etale-locally on $U$, there is an invertible $\Mc$ on $\Yc\times_S U$ such that $m=[\Mc]$ and $\Lc\simeq f^*\Mc$. This proves that $\psi$ is an epimorphism.

Saying that $\psi$ is a monomorphism just means that for any $U$, $\psi(U)$ is fully faithful. Let $\Mc_1$ and $\Mc_2$ be two invertible sheaves on $\Yc\times_S U$. We have to prove that
$$\Isom(\Mc_1, \Mc_2) \flechelongue \Hom(([\Mc_1], f^*\Mc_1),([\Mc_2],f^*\Mc_2))$$ 
is bijective. This is obvious if the right-hand side is empty. Otherwise we have $[\Mc_1]=[\Mc_2]$ in $\Pic_{\Yc/S}(U)$ and there is an isomorphism $f^*\Mc_1\simeq f^*\Mc_2$. Using~\cite{Brochard_Picard}~2.2.6, this implies that $\Mc_1$ and $\Mc_2$ are isomorphic. But then both sides of the above map naturally identify to $\gm(U)$, which concludes the proof.
\end{demo}

\begin{sousthm}\label{thm_ppal}
Let $S$ be an integral base scheme and let $\X$, $\Y$ be proper and finitely presented algebraic stacks over $S$. Let $f : \X \fleche \Y$ be a surjective morphism. Then there is a nonempty open subset $U$ of $S$ with the following properties:
\begin{itemize}
 \item[a)] The functors $(\pic)_{|_U}$ and $(\Pic_{\Y/S})_{|_U}$ are algebraic spaces, and the morphism $f^* : (\Pic_{\Y/S})_{|_U} \fleche (\pic)_{|_U}$ is quasi-affine and of finite type.
\item[b)] The stacks $\champic(\X/S)_{|_U}$ and $\champic(\Y/S)_{|_U}$ are algebraic stacks, and the morphism $f^* : \champic(\Y/S)_{|_U} \fleche \champic(\X/S)_{|_U}$ is of finite type with affine diagonal.
\end{itemize}
\end{sousthm}
\begin{sousremarque} \rm
\label{remarque_hyp_suppl}
If moreover $\X$ and $\Y$ are reduced with geometrically reduced and geometrically connected fibers, then the morphism $f^*$ in~b) is quasi-affine as well (lemma~\ref{comparaison}).
\end{sousremarque}
\begin{sousremarque}\rm We will see later (\ref{cor_sur_un_corps}) that if $s$ is a point of $S$, then the fiber morphism $f^*_s : (\Pic_{\Y/S})_s \fleche (\Pic_{\X/S})_s$ is actually affine. 
\end{sousremarque}
\begin{demo} We can assume that $S$ is affine, of finite type over $\Z$, and regular. Let us start with some further reductions.
\medskip

\begin{etape}{We can assume that $\X_{\red}$ and $\Y_{\red}$ have geometrically reduced and geometrically connected fibers.}
 Indeed, using the lemma~\ref{fgregc} and finite flat descent, we can assume that the connected components of $\X$ and $\Y$ have this property. Let us denote by $\X_i$ (resp. $\Y_j$) the connected components of $\X$ (resp. $\Y$). For every $i$, there is a $j$ such that $f(\X_i)\subset \Y_j$. Since $\pic=\prod_i \Pic_{\X_i/S}$ and $\Pic_{\Y/S}=\prod_j \Pic_{\Y_j/S}$ (and similarly for the Picard stacks) we can replace $\X$ by $\X_i$ and $\Y$ by $\Y_j$.
\end{etape}

\begin{etape}{We can assume that $\X$ and $\Y$ are reduced.}
Indeed, there is a commutative diagram
$$\xymatrix{
\Pic_{\Y/S} \ar[r]^{f^*}\ar[d] & \pic. \ar[d]\\
\Pic_{\Y_{\textrm{red}}/S} \ar[r]_{f_{\textrm{red}}^*} &\Pic_{\X_{\textrm{red}}/S}
 }$$
Shrinking $S$ if necessary, the vertical maps are affine and finite type morphisms of algebraic spaces owing to~\ref{cas_nilpotent}. Now, if $f_{\textrm{red}}^*$ is quasi-affine and of finite presentation, then so is $f^*$. The same diagram with Picard stacks instead of Picard functors shows that if $f^*_{\red}$ is of finite type with affine diagonal, then the same holds for $f^*$. From now on, to prove a) and b) we can work under the additional assumptions of~\ref{remarque_hyp_suppl}.
\end{etape}

\begin{etape}{We can assume that $P(f)$ holds (see the definition of $P(f)$ on p.\pageref{pagepdef}).}
 Indeed, let us consider the Stein factorization of $f$, 
$$\xymatrix{f : \X \ar[r]^-{f_0} & \Specc(f_*\Oc_{\X})=\Y_0 \ar[r]^-h & \Y\, .}$$
Then the morphism $\Oc_{\Y_0} \fleche (f_0)_*\Oc_{\X}$ is an isomorphism and thus $P(f_0)$ holds. Moreover, $h$ is finite, and since it is surjective and $\Y$ is reduced, we can apply lemma~\ref{lem_factorisation} and $h$ is the composition of a finite number of morphisms satisfying the property $P$. Note that $\Y_0$ and all the stacks given by~\ref{lem_factorisation} are reduced. Note also that, since $\X$ has geometrically connected fibers, and since all the morphisms in this factorization are surjective, it follows that all the stacks $\Y_i$ have geometrically connected fibers. Below, we will prove that $P(f)$ holds universally over a nonempty open subset of $S$. Since $\X$ has geometrically reduced fibers, this will imply that the fibers (over this open subset) of the stacks occurring in the above factorization are still geometrically reduced.
\end{etape}

\begin{etape}{We can assume that $P(f)$ holds universally.}
Let $A=\Oc_S$ and let $k$ be its global dimension. First, we notice that $P(f)$ remains true after any flat base change (\cite{Brochard_Picard}~A.3.4). Thus using generic flatness, we can assume that $\X$, $\Y$ and all the $R^if_*\Oc_{\X}$ and $R^ig_*\Oc_{\X\times_{\Y}\X}$ for $i\leq k$ (where $g$ is the natural morphism from~$\X\times_{\Y}\X$ to $\Y$) are flat over $S$. Then forming $f_*\Oc_{\X}$ and $g_*\Oc_{\X\times_{\Y}\X}$ commutes with any base change. Indeed, this question is local on $\Y$ so to prove it we can assume that $\Y$ is an affine scheme $\Spec B$. Now $g_*\Oc_{\X\times_{\Y}\X}$ and $f_*\Oc_{\X}$ are the coherent $\Oc_{\Y}$-modules corresponding to the $B$-modules $H^0(\X\times_{\Y}\X,\Oc_{\X\times_{\Y}\X})$ and $H^0(\X,\Oc_{\X})$, so we only have to prove that forming these $B$-modules commutes with any base change $A\fleche A'$, which follows immediately from~\ref{kunneth_tame}. Let us denote by $\pi_1$ and $\pi_2$ the projection maps from $\X\times_{\Y}\X$ to $\X$. By generic flatness again, we can assume that $\Coker(\pi_1^*-\pi_2^*)$ is flat over $S$. Now, we claim that $P(f)$ holds universally. Indeed, forming $\Coker(\pi_1^*-\pi_2^*)$ commutes with (and thus it remains flat after) any base change. This implies that the same holds for $\Im(\pi_1^*-\pi_2^*)$, hence also for $\Ker(\pi_1^*-\pi_2^*)$. Thus the latter is universally equal to $\Oc_{\Y}$.
\end{etape}

\begin{souslem}[\cite{SGA6} XII 4.2]
\label{lem_descente}
 Let $f : \X \fleche \Y$ be a morphism of algebraic stacks. Assume that
$P(f)$ holds and denote by $g$ the canonical map from $\X\times_{\Y}\X$ to $\Y$. Let $\Lc$, $\Mc$ be locally free sheaves of finite rank over $\Y$. Then the diagram
$$\xymatrix{\Hom_{\Oc_{\Y}}(\Lc,\Mc) \ar[r]& \Hom_{\Oc_{\X}}(f^*\!\Lc,f^*\!\Mc) \ar@<0.5ex>[r] \ar@<-0.5ex>[r]&
\Hom_{\Oc_{\X\times_{\Y}\X}}(g^*\!\Lc,g^*\!\Mc)}$$
is exact.
\end{souslem}
\begin{demo}
 This is an easy fact and the (short) proof given in~\cite{SGA6}~XII~4.2 works without any change.
\end{demo}

Let us now finish the proof of Theorem~\ref{thm_ppal}. First, we can assume that $\pic$ and $\Pic_{\Y/S}$ (resp. $\champic(\X/S)$ and $\champic(\Y/S)$) are both algebraic spaces (resp. algebraic stacks) locally of finite presentation over $S$ (use~\ref{champicqsep}~(1) and~\ref{lem_pic}). Moreover we can assume that the diagonal of~$\champic(\Xc\times_{\Yc}\Xc)$ is quasi-affine~(\ref{champicqsep}~(5)). Using the lemma~\ref{comparaison} and faithfully flat descent, it is enough to prove that the map~$f^*$ in b) is quasi-affine and of finite type. For this we have to prove that given an invertible sheaf~$\Lc$ on~$\Xc$, corresponding to a map~$S \fleche \champic(\Xc/S)$, the fiber product~$\Zc:=\champic(\Y/S)\times_{\champic(\X/S)} S$ is a quasi-affine scheme of finite presentation (assuming that $S$ is affine and noetherian, but not integral nor regular anymore).


The stack $\Zc$ can be described as follows. An object of the fiber category $\Zc(T)$ is a couple $(\Mc,\alpha)$ where $\Mc$ is an invertible sheaf on $\Yc\times_S T$ and $\alpha$ is an isomorphism $f^*\Mc \fleche \Lc_T$. Note that because of~\ref{comparaison}~$\Zc$ is actually an algebraic space.

Let us denote by $I$ the sheaf $\fIsom(\pi_1^*\Lc,\pi_2^*\Lc)$ on $S$-schemes. Because of our previous reductions on the diagonal of~$\champic(\Xc\times_{\Yc}\Xc)$, $I$ is a quasi-affine scheme of finite presentation over~$S$.
Given an invertible sheaf $\Mc$ on $\Y$ and an isomorphism $\alpha : f^*\Mc\fleche \Lc$, the canonical isomorphism $\pi_1^*f^*\!\Mc \fleche \pi_2^*f^*\!\Mc$ induces an isomorphism from $\pi_1^*\Lc$ to $\pi_2^*\Lc$, thus an element of $I(S)$. This construction is clearly functorial in $S$ and yields a morphism $\Zc \fleche I$.

Let us prove that this morphism is a monomorphism. Consider two points $(\Mc_1,\alpha_1)$ and  $(\Mc_2,\alpha_2)$ which have the same image in $I(S)$. Then we have an isomorphism $\alpha_2^{-1}\alpha_1$ from $f^*\Mc_1$ to $f^*\Mc_2$. Saying that $(\Mc_1,\alpha_1)$ and  $(\Mc_2,\alpha_2)$ have the same image in $I(S)$ means that in the diagram
$$\xymatrix@C=1.7pc{\Hom_{\Oc_{\Y}}\!(\Mc_1,\Mc_2) \ar[r]& \Hom_{\Oc_{\X}}\!(f^*\!\!\Mc_1,f^*\!\!\Mc_2) \ar@<0.5ex>[r]^-{\pi_1^*} \ar@<-0.5ex>[r]_-{\pi_2^*}& \Hom_{\Oc_{\X\!\times_{\Y}\!\X}}\!(g^*\!\!\Mc_1,g^*\!\!\Mc_2)}$$
we have ${\pi_1^*}(\alpha_2^{-1}\alpha_1)={\pi_2^*}(\alpha_2^{-1}\alpha_1)$. Using the exactness~(\ref{lem_descente}) of this diagram, and applying the same argument with $\alpha_1^{-1}\alpha_2$, we deduce that there is an isomorphism $\Mc_1 \fleche \Mc_2$ over $\alpha_2^{-1}\alpha_1$. This proves that $(\Mc_1,\alpha_1)$ and  $(\Mc_2,\alpha_2)$ are isomorphic in $\Zc(S)$, as desired.

Now, applying the lemma~\ref{lem_quaff}, we see that $\Zc\fleche I$ is quasi-affine and of finite presentation, and this concludes the proof.
\end{demo}

To conclude this section, let us state two immediate corollaries.

\begin{souscor}\label{cor_repres}
Let $\X$ be a proper and finitely presented algebraic stack over an integral base scheme $S$. Then there is a nonempty open subset $V$ of $S$ such that the Picard functor $(\pic)_{|_V}$ is representable by a scheme that is a disjoint union of quasi-projective and finite presentation open subschemes.
\end{souscor}
\begin{demo}
Using~\cite{Olsson_Hom_stacks}~2.2 we can assume that $S$ is noetherian. Now the result is known if $\X$ is a scheme (see \cite{SGA6} XII 1.2). In the general case, let $X\fleche \X$ be a Chow presentation of $\X$, \emph{i.e.} a proper and surjective morphism where $X$ is a scheme. Then replacing $S$ with a nonempty open subscheme, we may assume that $\Pic_{X/S}$ is representable by a scheme which is a disjoint union of quasi-projective and finite presentation open subschemes (\emph{loc. cit.}), and that the morphism $\pic \fleche \Pic_{X/S}$ is representable by a quasi-affine and finite presentation morphism~(\ref{thm_ppal}). This gives the result.
\end{demo}

\begin{souscor}\label{cor_sur_un_corps} $ $
\begin{itemize}
\item[(i)] Let $S$ be the spectrum of a field and $\X$ a proper algebraic stack over $S$. Then the relative Picard functor $\pic$ is representable by a scheme which is a disjoint union of quasi-projective and finite presentation open subschemes.
\item[(ii)] Let $f : \X \fleche \Y$ be a surjective morphism between two proper algebraic stacks over a field $k$. Then the induced morphism $$f^* : \Pic_{\Y/k} \flechelongue \Pic_{\X/k}$$ is affine.
\end{itemize}
\end{souscor}
\begin{demo} (i) follows immediately from~\ref{cor_repres}. For~(ii), we know by~\ref{thm_ppal} that $f^*$ is quasi-affine. Using~\cite{SGA6}~XII~1.4, we see that $f^*$ is actually affine.
\end{demo}

\section{Corollaries, examples and applications}
\label{Corollaries}

To state the results concerning the Picard functors, we need some finiteness notions for non-representable functors or non-algebraic stacks. The paragraph~\ref{prelim} is devoted to the definitions and elementary properties of these notions. 

\subsection{Finiteness properties for non-representable functors}
\label{prelim}
Let $S$ be a base scheme. We consider functors from $(Sch/S)^{\text{op}}$ to $(Sets)$. We will need the notions of quasi-compact, quasi-separated or finite presentation functor or morphim of functors, even if the functors are not representable.

\begin{sousdefi}
\label{def_surj}
A morphism of functors $F \fleche G$ is surjective if for every field $K$ and every point $\xi$ in $G(K)$, there is an extension $L$ of $K$ such that $\xi_L$ is in the image of $F(L)$. \end{sousdefi}

\begin{sousdefi} $ $
\label{def_qcqsep}
\begin{itemize}
\item[(i)] A functor $F$ is said to be quasi-compact if there exists a quasi-compact $S$-scheme $T$ and a surjective morphism $T\fleche F$.
\item[(ii)] A morphism of functors $F \fleche G$ is said to be quasi-compact if for every quasi-compact scheme $T$ and every morphism $T\fleche G$, the functor $F_T$ obtained by base change is quasi-compact.
\item[(iii)] A morphism of functors $F \fleche G$ is said to be quasi-separated if the diagonal morphism $F\fleche F\times_G F$ is quasi-compact.
\item[(iv)] 
A functor $F$ is said to be locally of finite presentation (implicitly, over~$S$) if for every filtering inverse system $(Z_{\lambda})_{\lambda\in\Lambda}$ of $S$-schemes, such that each $Z_{\lambda}$ is an affine scheme, the map
$$\lind F(Z_{\lambda}) \flechelongue F(\lpro Z_{\lambda})$$
is bijective.
\item[(v)] A morphism of functors $F\fleche G$ is said to be locally of finite presentation if for every $S$-scheme $U$ and every morphism $U\fleche G$, the functor $F_U=F\times _G U$ is locally of finite presentation over $U$.
\item[(vi)] A morphism of functors $F\fleche G$ is said to be of finite presentation if it is locally of finite presentation, quasi-compact, and quasi-separated.
\end{itemize}
\end{sousdefi}

The following properties are straightforward. We give detailed proofs in the note~\cite{Brochard_qcqsep}.

\begin{sousprop} \label{pptes_qcqsep} $ $
\begin{itemize}
\item[(i)] If $F$ is an algebraic space, or if the morphism $F\fleche G$ is representable by algebraic spaces, the above notions coincide with the usual ones.
\item[(ii)] Every isomorphism is quasi-compact.
\item[(iii)] Every monomorphism is quasi-separated.
\item[(iv)] The class of surjective (resp. quasi-compact, quasi-separated, finite presentation, locally of finite presentation) morphisms is stable by base change and composition.
\item[(v)] If $F\fleche G$ is a surjective morphism and if $F$ is quasi-compact, then $G$ is quasi-compact.
\item[(vi)] If $F\fleche G$ is a quasi-compact morphism and if $G$ is quasi-compact, then $F$ is also quasi-compact.
\item[(vii)] Let $$\xymatrix@C=1pc@R=1pc{F \ar[rr]^f \ar[rd]_h && G\ar[ld]^g\\ &H}$$ be a commutative diagram of functors. If $h$ is quasi-compact and $g$ is quasi-separated, then $f$ is quasi-compact.
\item[(viii)] In a diagram as in (vii), if $h$ is quasi-separated, then $f$ is too.
\item[(ix)] Let $$\xymatrix{F' \cartesien \ar[r]^{\varphi'} \ar[d]_{f'} &F\ar[d]^f\\ G'\ar[r]^{\varphi}&G}$$ be a cartesian diagram of functors. Assume that the base change morphism $\varphi$ is surjective and quasi-compact. Then $f$ is quasi-compact (resp. quasi-separated), if and only if $f'$ is.
\item[(x)] In a diagram as in (vii), if $g$ and $h$ are locally of finite presentation, then $f$ is too. If moreover $h$ is of finite presentation and $g$ is quasi-separated, then $f$ is of finite presentation.
\end{itemize}
\end{sousprop}

\begin{sousremarque}\rm \label{comp_typefini_raynaud}
We explain in the note \cite{Brochard_qcqsep} that the definition of a quasi-compact morphism given above generalizes Kleiman's definition of a finite type morphism of Picard functors (SGA6~\cite{SGA6}~XIII).
\end{sousremarque}

\begin{sousremarque}\rm \label{qcqsep_champs}
For a (non-necessarily algebraic) stack over $S$, we have analogous notions to~\ref{def_surj} and~\ref{def_qcqsep}, with the same formal properties as in~\ref{pptes_qcqsep}. We leave the details to the reader\footnote{There is however a little problem of terminology. Usually, an algebraic stack is said to be quasi-separated if its diagonal is quasi-compact and separated. With the definitions above the diagonal is only quasi-compact (note that this problem did not occur with algebraic spaces since the diagonal of an algebraic space is always separated). To avoid confusion in the sequel, we will stick to the usual definition and only talk about ``stacks with quasi-compact diagonal'' when we need to.}.
\end{sousremarque}

\subsection{General finiteness results for the Picard functor}
\label{general}

The propositions \ref{pic_qsep} and \ref{morph_qcqsep} below are easy consequences of~\ref{thm_ppal}, using noetherian induction. Note that~\ref{morph_qcqsep}~(ii) and~\ref{puiss_n} generalize~\cite{SGA6}~XIII, (3.5) and~(3.6).

\begin{sousprop} \label{pic_qsep}
Let $\X$ be a proper and finitely presented algebraic stack over a base scheme $S$. Then the Picard functor $\pic$ is quasi-separated, and the Picard stack $\champic(\X/S)$ has a quasi-compact diagonal (in the sense of~\ref{def_qcqsep} and~\ref{qcqsep_champs}).
\end{sousprop}
\begin{demo}
By standard limit arguments, we can assume that $S$ is noetherian (\cite{Olsson_Hom_stacks}~2.2). Let $T$ denote the disjoint union of the irreducible components of $S$. The morphism $T\fleche S$ is surjective and quasi-compact, so that using~\ref{pptes_qcqsep}~(ix) we may assume $S$ is irreducible. Plainly we can assume $S$ is integral. By noetherian induction it is enough  to prove that there is a nonempty open subset $U$ of $S$ such that $\pic \times_S U$ is quasi-separated (apply~\ref{pptes_qcqsep}~(ix) with the base change $U\amalg (S\setminus U) \fleche S$). This fact is an immediate consequence of~\ref{cor_repres}. The same proof works for the Picard stack.
\end{demo}

\begin{sousprop}\label{morph_qcqsep}$ $
\begin{itemize}
\item[(i)] Let $\X$ and $\Y$ be algebraic stacks over a scheme $S$. Assume that $\Y\!$ is proper and finitely presented. Then every $S$-morphism $\Pic_{\Y\!/S} \fleche \Pic_{\X\!/S}$ is quasi-separated and locally of finite presentation. Moreover, every $S$-morphism from $\champic(\Y/S)$ to~$\champic(\X/S)$ has a quasi-compact diagonal and is locally of finite presentation.
\item[(ii)] Let $\X$ and $\Y$ be proper and finitely presented algebraic stacks over a base scheme $S$. Let $f : \X \fleche \Y$ be a surjective morphism. Then the morphisms $f^* : \Pic_{\Y\!/S} \fleche \pic$ and $f^* : \champic(\Y\!/S) \fleche \champic(\X\!/S)$ are of finite presentation.
\end{itemize}
\end{sousprop}
\begin{demo}
(i) Apply \ref{pic_qsep} and \ref{pptes_qcqsep} (viii). The finite presentation assertion is well known.

(ii) We already know that $f^*$ is quasi-separated and locally of finite presentation. It remains to prove that it is quasi-compact. As in the proof of~\ref{pic_qsep}, using~\ref{pptes_qcqsep}~(ix) we can assume that $S$ is noetherian and integral. By noetherian induction it is enough to prove the assertion over a nonempty open subset of $S$. It is now an immediate consequence of~\ref{thm_ppal}. The same proof works for $f^* : \champic(\Y\!/S) \fleche \champic(\X\!/S)$.
\end{demo}

\begin{sousprop}\label{puiss_n}
Let $\X$ be a proper and finitely presented algebraic stack over a scheme $S$. Then for every positive integer $n$, the morphisms
$$\fonction{\varphi_n}{\pic}{\pic}{\Lc}{\Lc^{\otimes n}} \quad \text{and} \quad
\fonction{\lambda_n}{\champic(\X/S)}{\champic(\X/S)}{\Lc}{\Lc^{\otimes n}}$$
are of finite presentation.
\end{sousprop}
\begin{demo}
Again we only have to prove that $\varphi_n$ and $\lambda_n$ are quasi-compact. By Chow's lemma (\cite{Olsson_lemme_chow}~1.1) there is a proper and surjective morphism $\pi : X \fleche \X$ with $X$ a projective $S$-scheme. We have a commutative diagram:
$$\xymatrix{\pic \ar[r]^{\pi^*} \ar[d]_{\varphi_{n,\X}} & \Pic_{X/S} \ar[d]^{\varphi_{n,X}}\\
\pic \ar[r]_{\pi^*} & \Pic_{X/S}.}$$
The morphism $\pi^*$ is quasi-compact by~\ref{morph_qcqsep}, and so is~$\varphi_{n,X}$ by~\cite{SGA6}~XIII~3.6. So $\varphi_{n,\X}$ is quasi-compact by~\ref{pptes_qcqsep}~(vii).

To prove that $\lambda_n$ is quasi-compact, as in the preceding proofs we can assume that $S$ is integral, and by noetherian induction it is enough to prove that the assertion holds over a nonempty open subset of $S$. We can obviously assume that $\X\fleche S$ is surjective and that $\X$ is connected. Using~\ref{cas_nilpotent} and the diagram
$$\xymatrix{\champic(\X/S) \ar[r]^{\lambda_n} \ar[d] & \champic(\X/S)\ar[d] \\
 \champic(\X_{\red}/S)\ar[r]_{\lambda^{\red}_n}&  \champic(\X_{\red}/S)}$$
we can assume that $\X$ is reduced. With these remarks and with~\ref{fgregc}, we can restrict to the case where $\X$ has geometrically reduced and geometrically connected fibers. By \emph{fppf} descent we can also assume that $\X$ has a section. Now, by~\cite{Brochard_Picard}~2.3.4, the stack $\champic(\X/S)$ is isomorphic to $\pic \times_S \bgm$. Hence $\lambda_n$ is the product of $\varphi_n$ and the $n^{\textrm{th}}$-power map $\bgm \fleche \bgm\,$, so it suffices to prove that the latter is quasi-compact, which is obvious since $\bgm$ itself is quasi-compact.
\end{demo}

\subsection{The torsion component $\pictox$}
\label{torsion}

We recall the definition of the torsion and the neutral component of a (commutative) group functor.

\begin{sousdefi}
Let $G : (Sch/S)^{\text{op}} \fleche (Gr)$ be a functor, where $(Gr)$ is the category of commutative groups. We define two group subfunctors $G^0$ and $G^{\tau}$ of $G$ as follows. If $T$ is the spectrum of an algebraically closed field $K$, we say that a point $t\in G(T)$ is in $G^0(T)$ if there are connected $K$-schemes of finite type $T_1, \dots, T_n$, and for each $i$ a morphism $\alpha_i : T_i\fleche G$ (\emph{i.e.} $\alpha_i\in G(T_i)$) and two $K$-points $s_i, t_i \in T_i$ such that
\beqn \alpha_1(s_1) &=&t\\
\alpha_1(t_1)&=& \alpha_2(s_2) \\
\vdots &&\\
\alpha_{n-1}(t_{n-1})&=&\alpha_n(s_n)\\
\alpha_n(t_n)&=&0.
\eeqn
For an arbitrary $S$-scheme $T$, we say that a point $t\in G(T)$ is in $G^0(T)$ (resp. in $G^{\tau}(T)$) if for every algebraically closed field $K$ and every morphism $\xi : \Spec K \fleche T$, the restriction $t_{\xi}\in G(\Spec K)$ is in $G^0(\Spec K)$ (resp. if there is $n>0$ such that $t^n\in G^0(T)$).
\end{sousdefi}

We let the reader check that these are indeed group subfunctors of $G$, and that forming them commutes with any base change. Moreover, if $G$ is a group scheme locally of finite type over a field $k$, then $G^0$ is the connected component of the neutral element (hence it is open, closed, geometrically irreducible and of finite type). In this case $G^{\tau}$ is then an open and closed group subscheme of $G$: it is open because it is the union $\cup_{n>0} \varphi_n^{-1}(G^0)$ where $\varphi_n : G\fleche G$ is the multiplication by $n$. To prove that it is closed, by faithfully flat descent we may assume that the base field is algebraically closed, and then the complement of $G^{\tau}$ is a union of translates of $G^{\tau}$, more precisely
$$ G\setminus G^{\tau}=
\bigcup_{g\in G(k)\setminus G^{\tau}(k)}gG^{\tau}\, .$$
 The reader can also check that $G^{\tau}$ contains any group subscheme of finite type of $G$.

For an arbitrary base scheme $S$, $G^0$ may not be representable (but its fibers are, by the above). In the case of the Picard functor $\pic$ of a proper algebraic stack $\X$ over a base scheme $S$, we have proved in~\cite{Brochard_Picard} (4.2.10) that the subfunctor $\piczero$ is representable if $\pic$ is a locally finitely presented algebraic space that is smooth along the unit section. In this case $\piczero$ is an open subspace of $\pic$, and is of finite type over $S$. We prove in~\ref{picto_ouvert} below that the same conclusions hold for $\pictox$ without the additional assumption on $\pic$.

\begin{souslem}\label{image_inverse_picto}
Let $\X$ and $\Y$ be proper algebraic stacks over a noetherian scheme $S$, and $f : \X \fleche \Y$ a surjective morphism. Then the natural diagram
$$\xymatrix{\pictoy \ar[r] \ar[d] & \pictox \ar[d] \\
\Pic_{\Y/S} \ar[r]_{f^*} & \pic}$$
is cartesian. In other words, for every $S$-scheme $T$ and every element $l$ in $\Pic_{\Y/S}(T)$, $l$ is in $\pictoy(T)$ if and only if $f^*(l)$ is in $\pictox(T)$.
\end{souslem}
\begin{demo}
For any geometric point $\Spec K \fleche S$ the functors $(\Pic_{\Y\!/S})_K$ and $(\Pic_{\X\!/S})_K$ are schemes locally of finite type~(\ref{cor_sur_un_corps}) over $K$ and the morphism $$f^* : (\Pic_{\Y/S})_K \fleche (\pic)_K$$ is quasi-compact~(\ref{morph_qcqsep}). Thus we can apply~\cite{SGA6}~XIII~4.2.
\end{demo}

\begin{sousthm}\label{picto_ouvert}
Let $\X$ be a proper and finitely presented algebraic stack over a base scheme $S$. Then:
\begin{itemize}
\item[(i)] The morphism $\pictox \fleche \pic$ is representable by an open immersion.
\item[(ii)] $\pictox$ is of finite presentation over $S$.
\end{itemize}
\end{sousthm}
\begin{demo}
We can assume that $S$ is noetherian. By Chow's lemma, there is a projective scheme $X$ over $S$ and a surjective morphism $\pi : X \fleche \X$. By~\ref{image_inverse_picto}, we have a cartesian square
$$\xymatrix{\pictox \ar[r]^{\pi^{*,\tau}} \ar[d]_{i_{\X}} \cartesien & \picto_{X/S} \ar[d]_{i_{X}} \\
\pic \ar[r]^{\pi^{*}} & \Pic_{X/S}.}$$
(i) Applying \cite{SGA6}~XIII, (4.7) (i) to the scheme $X/S$, we see that the morphism $i_X$ is representable by an open immersion. So $i_{\X}$ is too.

(ii) By \ref{morph_qcqsep} the morphism $\pi^*$ is quasi-compact, hence $\pi^{*,\tau}$ is quasi-compact too. By \cite{SGA6}~XIII~4.7 (iii), $\picto_{X/S}$ is quasi-compact over $S$. So $\pictox$ is quasi-compact over~$S$. Moreover $i_{\X}$ is quasi-separated and locally of finite presentation since it is an open immersion. But $\pic$ is quasi-separated and locally of finite presentation over~$S$. Then so is $\pictox$. \end{demo}

\begin{sousexemple} \rm
\label{exemple_piczero_non_repres} 
The latter theorem provides an open subgroup of $\pic$ which is of finite type over $S$. This can be very useful when the existence theorems for the neutral component $\piczero$ do not apply. Let us give an example in which $\piczero$ is not representable. Let $S=\Spec \Z$ and let $\X$ be the classifying stack $B(\Z/n\Z)_{\P^1_{\Z}}$ over $\P^1_{\Z}$. In other words, $\X=B(\Z/n\Z)\times_S \P^1_{\Z}$. Note that the stack $\X$ is smooth and proper over $S$, with geometrically integral fibers. If $T$ is a scheme over $S$, we have an isomorphism (see~\cite{Brochard_Picard}~Exemple~5.3.7)
$$\Pic(\X\times_S T) \simeq \Pic(\P^1_T)\times \widehat{(\Z/n\Z)}(\P^1_T)$$
where $\widehat{(\Z/n\Z)}=\Homs_{S-\text{Gr}}(\Z/n\Z,\gm)$ is the Cartier dual of $\Z/n\Z$.
But $\widehat{(\Z/n\Z)}\simeq \mu_n$ and $\mu_n(\P^1_T)\simeq \mu_n(T)$, hence the Picard functor of $\X/S$ is representable and we have an isomorphism
$$\pic \simeq \Z\times \mu_n.$$
As a scheme, $\pic$ is a disjoint union of copies of $\mu_n$, indexed by $\Z$ (in particular it is not of finite type). The open subgroup $\pictox$ provided by~\ref{picto_ouvert} is the $0^{\text{th}}$ copy of~$\mu_n$. It is a finite flat group scheme over $S$. On the other hand, the subfunctor $\piczero$ coincides with the neutral component $\mu_n^0$ of $\mu_n$. The reader can check easily that it is not an open subspace of $\mu_n$.

\end{sousexemple}

\begin{sousexemple}\rm
 In the case of root stacks, we can also prove directly Theorem~\ref{picto_ouvert} from~\ref{puiss_n}. Let $X$ be an $S$-scheme that is proper (resp. projective), flat and with geometrically integral fibers. Let $\Lc$ be an invertible sheaf on $X$. We know that $\Pic_{X/S}$ is an algebraic space. Let us denote by~$\X=[\Lc^{\frac1n}]$ the stack, the fiber category over an $S$-scheme $U$ of which is the category of triples $(x,\Mc,\varphi)$ where $x\in X(U)$, $\Mc$ is an invertible sheaf on $X$ and $\varphi$ is an isomorphism of invertible sheaves from $\Mc^{\otimes n}$ to $x^*\Lc$. Then by~\cite{Brochard_Picard}~5.3, there is a short exact sequence of \'etale sheaves:
$$\xymatrix{0 \ar[r] &\Pic_{X/S} \ar[r]& \pic \ar[r] &\Z/n\Z \ar[r] & 0.}$$
Moreover, $\pic$ is an algebraic space, and $\pi^* : \Pic_{X/S} \fleche \pic$ is an open and closed immersion. As in~\ref{puiss_n}, let us denote by $\varphi_n$ the  $n^{\textrm{th}}$-power map of the Picard functor $\pic$. Because of the above exact sequence, the image of $\varphi_n$ lies in~$\Pic_{X/S}$. We deduce easily from this remark, and from the fact that $\pi^*$ induces an isomorphism from $\Pic^0_{X/S}$ to $\piczero$ (see \emph{loc. cit.}), that there is a cartesian square:
$$\xymatrix{\pictox \cartesien \ar[r] \ar[d]& \picto_{X/S} \ar[d]^g \\
\pic \ar[r]_{\varphi_n} & \pic}$$
where $g$ is the composition of $\pi^*$ and of $i : \picto_{X/S} \fleche \Pic_{X/S}$. We know by~\cite{SGA6}~XIII~4.7 that $i$ is an open (resp. open and closed) immersion and that $\picto_{X/S}$ is of finite type over $S$. It follows that $\pictox \fleche \pic$ is an open (resp. open and closed) immersion. Since $\varphi_n$ is of finite type by~\ref{puiss_n}, we also deduce that $\pictox$ is of finite type over $S$ as well.

For a concrete example, if $X=\P_S^k$ and $\Lc=\Oc(l)$, then $\Pic_{X/S}^0$, $\picto_{X/S}$ and $\piczero$ are trivial $S$-groups. Let $d=\gcd(l,n)$. We can see easily that $\pictox$ is isomorphic to the constant group scheme $\Z/d\Z$, generated by the section corresponding to $\Omega^{\frac{n}d}\otimes \Oc(-\frac{l}d)$, where $\Omega$ is the canonical $n^{\textrm{th}}$ root of $\Lc$ on $\X$.
\end{sousexemple}


%
%
%

\subsection{Arithmetic results}
\label{arithmetic}

\begin{sousthm}\label{Neron_Severi}
Let $\X$ be a proper algebraic stack over a noetherian base scheme $S$. If $\xi : \Spec K \fleche S$ is a geometric point of $S$, let us denote by $\nersev(\xi)$ (or $\nersev(\X, \xi)$ if we need to precise the stack) the N\'eron-Severi group of the geometric fiber $\X_{\xi}$ of $\X$ over $\xi$.
$$\nersev(\xi) = \frac{\Pic_{\X_{\xi}/K}(K)}{\Pic^0_{\X_{\xi}/K}(K)}\, .$$
Then the groups $\nersev(\xi)$ are of finite type. Moreover, their rank and the order of their torsion subgroups are uniformly bounded over $S$.
\end{sousthm}
\begin{demo} It is enough to prove the second assertion. To bound the rank, let us take a Chow presentation $\pi : X\fleche \X$. By~\ref{image_inverse_picto}, we have a cartesian square of group functors
$$\xymatrix{\disp \pictox \ar[r] \ar[d] \cartesien & \picto_{X/S} \ar[d] \\
\pic \ar[r] & \Pic_{X/S}.}$$ 
This implies that for each geometric point $\xi$, the morphism $\pi^*$ induces an injective morphism
$$\xymatrix{\disp \frac{\nersev(\X, \xi)}{\nersev(\X, \xi)_{\rm tors}} = \frac{\Pic_{\X_{\xi}/K}(K)}{\picto_{\X_{\xi}/K}(K)} \ar@{^(->}[r] & \disp
\frac{\Pic_{X_{\xi}/K}(K)}{\picto_{X_{\xi}/K}(K)}
=\frac{\nersev(X, \xi)}{\nersev(X, \xi)_{\rm tors}}} .$$
But the group on the right is free of finite rank uniformly bounded over $S$ (\cite{SGA6}~XIII~5.1). Hence the same holds for the left-sided group (with the same bound as for $X$).

To see that the order of the torsion subgroup is bounded, the proof given in \emph{loc.cit.} works without any change: using~\ref{cor_repres}, \ref{picto_ouvert} and noetherian induction, we may assume that $S$ is integral and that $\pictox$ is a scheme of finite type. For a geometric point $\xi$ of $S$, the torsion subgroup $\nersev(\xi)_{\textrm{tors}}$ is equal to the group $\picto_{\X_{\xi}/K}(K)/\Pic^0_{\X_{\xi}/K}(K)$, so its order is equal to the number of connected components of $\picto_{\X_{\xi}/K}$. But by~\cite{EGA} EGA~$\textrm{IV}_4$~9.7.9, there is a nonempty open subset of $S$ over which this number is constant.
\end{demo}

\begin{sousexemple}\rm
 Let $S$ be the spectrum of an algebraically closed field $k$, let $\Cc$ be a smooth $n$-pointed twisted curve in the sense of Abramovich and Vistoli, and let $C$ be its coarse moduli space. Let us denote by $(d_1,\dots,d_n)$ the $n$-uplet of integers defining the actions on the marked points, and by $\Lc_i$ the invertible sheaf $\Oc_{\Cc}(\Sigma_i)$ where $\Sigma_i \fleche \Cc$ is the $i^{\textrm{th}}$ marked point (see also~\cite{Cadman_USTITCOC}~2.2.4 and 4.1). Then by~\cite{Cadman_USTITCOC}, we know that there is a short exact sequence
$$\xymatrix{0 \ar[r] &\Pic(C) \ar[r]& \Pic(\Cc) \ar[r] & \disp \prod_{i=1}^n \Z/d_i\Z \ar[r]& 0,}$$
in which the class of the invertible sheaf $\Lc_i$ on $\Cc$ is sent to a generator of the factor $\Z/d_i\Z$.
Since $k$ is algebraically closed, we have $\Pic_{\Cc/S}(k)=\Pic(\Cc)$ and $\Pic_{C/S}(k)=\Pic(C)$. Moreover, by~\cite{Brochard_Picard} we know that $\Pic_{\Cc/S}^0=\Pic_{C/S}^0$. It follows that the above exact sequence induces an exact sequence involving the N\'eron-Severi groups:
$$\xymatrix{0 \ar[r] &\nersev(C) \ar[r]& \nersev(\Cc) \ar[r] & \disp \prod_{i=1}^n \Z/d_i\Z \ar[r]& 0.}$$
So $\nersev(\Cc)$ is the group obtained from $\nersev(C)$ by adding formally, for every $i$, a $d_i^{\textrm{th}}$ root of the class of $\Oc(D_i)$, where $D_i$ is the $i^{\textrm{th}}$ marked point on $C$. From the latter exact sequence, we  see that $\nersev(C)$ and $\nersev(\Cc)$ have the same rank. Moreover, the order of $\nersev(\Cc)_{\textrm{tors}}$ is bounded by $|\nersev(C)_{\textrm{tors}}|.\prod_{i=1}^n d_i$. More precisely, if $\iota$ is the index of the free part of $\nersev(C)$ in the free part of $\nersev(\Cc)$, we have the relation
$$\iota\times \frac{|\nersev(\Cc)_{\textrm{tors}}|}{|\nersev(C)_{\textrm{tors}}|}=\prod_{i=1}^n d_i\, .$$
\end{sousexemple}


In the rest of this section, we generalize to algebraic stacks some results of the article~\cite{Kahn_picfini} of Kahn.


\begin{souslem}[\cite{Kahn_picfini} th\'eor\`eme 1 a)]
Let $\X$ be a reduced algebraic stack, of finite type over $\Spec \Z$. Then the group $H^0(\X, \gm)$ is of finite type.
\end{souslem}
\begin{demo}
Let $X\fleche \X$ be a smooth and quasi-compact presentation of $\X$. Then $H^0(\X, \gm)$ is a subgroup of $H^0(X, \gm)$, which is of finite type owing to~\cite{Kahn_picfini} thm 1, a).
\end{demo}

\begin{souslem}
Let $\X \fleche \Y$ be a smooth, surjective and quasi-compact morphism of reduced algebraic stacks. Assume that $\X$ and $\Y$ are of finite type over $\Spec \Z$. If $\Pic(\X)$ is of finite type, then $\Pic(\Y)$ is of finite type.
\end{souslem}
\begin{demo} A descent argument (see for instance~\cite{Brochard_Picard}~2.1.2) shows that the kernel of the morphism $\Pic(\Y)\fleche \Pic(\X)$ is equal to the homology of the complex
$$\xymatrix@C=3.2pc{H^0(\X\!, \gm)\ar[r]^-{p_1^*-p_2^*} & H^0(\X\!\times_{\Y}\!\X\!, \gm) \ar[r]^-{p_{23}^*-p_{13}^*+p_{12}^*}& H^0(\X\!\times_{\Y}\! \X\!\times_{\Y}\! \X\!, \gm).}$$
The result follows since $H^0(\X\!\times_{\Y}\!\X\!, \gm)$ is of finite type.
\end{demo}

\begin{sousthm}[\emph{cf.} \cite{Kahn_picfini} th. 1, cor. 2 and cor. 4]\label{picfini}
$ $
\begin{itemize}
\item[a)] Let $S$ be equal to $\Spec \Z$ or $\Spec k$, where $k$ is a field of finite type over its prime subfield. If $\X$ is a normal and reduced algebraic stack of finite type over $S$, then the group $\Pic(\X)$ is of finite type.
\item[b)] In the case where $S$ is $\Spec k$, assume moreover that $\X$ is proper, geometrically normal and geometrically integral. Let $\ov{k}$ be an algebraic closure of $k$ and $G=\Gal(\ov{k}/k)$. Then the group $\Pic(\X_{\ov{k}})^G$ is of finite type (where $\X_{\ov{k}}$ denotes the fiber product $\X\times_{\Spec k} \Spec \ov{k}$).
\end{itemize}
\end{sousthm}
\begin{demo}
a) For the case $S=\Spec \Z$, consider a presentation $X\fleche \X$ of $\X$, with $X$ a normal scheme, of finite type over $S$. Then $\Pic(X)$ is of finite type because of the theorem 1 of \cite{Kahn_picfini}. The preceding lemma then implies that $\Pic(\X)$ is of finite type.

For the case $S=\Spec k$, the proof given in~\cite{Kahn_picfini} works without any change. We include it here for the convenience of the reader: let us choose an integral scheme $U$ of finite type over $\Spec \Z$, whose function field is $k$. Then by \cite{LMB}~4.18, replacing $U$ with a nonempty open subset if necessary, there is an algebraic stack $\Xt$ of finite presentation over $U$ whose generic fiber is equal to $\X$. Moreover we may assume $\Xt$ is normal. Now $\Pic(\Xt)$ is of finite type by the preceding case. Moreover, the morphism from $\Pic(\Xt)$ to $\Pic(\X)$ is surjective. This concludes the proof.

b) Let $P$ be the scheme $\Pic_{\X/k}$ (\ref{cor_sur_un_corps}) and $P^0$ its neutral component. Then $P^0$ is proper (\cite{Brochard_Picard}), thus an abelian variety. We have an exact sequence
$$\xymatrix{0\ar[r] & P^0(\ov{k}) \ar[r] & P(\ov{k}) \ar[r]& \nersev(\ov{k}),}$$
inducing
$$\xymatrix{0\ar[r] & P^0(\ov{k})^G \ar[r] & P(\ov{k})^G \ar[r]& \nersev(\ov{k}).}$$
Now, the last group is of finite type by~\ref{Neron_Severi}, and $P^0(\ov{k})^G$ is equal to $P(k)$, thus of finite type by~\cite{Kahn_picfini}~cor.~3. Moreover $P(\ov{k})=\Pic(\X_{\ov{k}})$.
\end{demo}

\sectionmark{Appendix}
\section{\!\!\! Appendix: Cohomology and base change for stacks}
\sectionmark{Appendix}
\label{appendice}

In scheme theory, the key point to get some ``base change theorems'' for cohomology of a coherent sheaf $\Fc$ on a proper scheme (over a noetherian ring~$A$), is the existence of a finite complex of finite free $A$-modules computing ``universally'' the cohomology of $\Fc$ (see~\cite{Mumford_Abelian_Varieties}~\S 5). A significant difference between an algebraic stack and a scheme (or even an algebraic space) is that the former can have infinite cohomological dimension. Consequently, there does not always exist such a finite complex. However, we explain below that even without such a complex, the semicontinuity theorem still holds for stacks (Theorem~\ref{semicont}).

Moreover, it is worth mentioning that there are two particular cases in which we can say much more:
\begin{itemize}
 \item[a)] if the stack is tame;
\item[b)] if the base ring has finite global cohomological dimension (\emph{e.g.} if it is regular and finite-dimensional).
\end{itemize}
Indeed, in these two cases, we can actually compute the cohomology of a given $\Oc_{\X}$-module $\Fc$ (at least the first cohomology groups) with a finite complex of $A$-modules. The existence of such complexes is provided by the Proposition~\ref{complexe_coh_tame} for the case~a), and by the Corollary~\ref{complexe_pour_base_reguliere} for~b).
Consequently, in these cases, the cohomology of $\Fc$ really behaves under base change as if $\X$ were a scheme (see~\ref{remarque_coh_tame}, \ref{cor_pour_base_reguliere}, \ref{other_csq}).

Since the proof of the semicontinuity theorem~\ref{semicont} relies on the case~b), we will start with these two particular cases.

\subsection{If the stack is tame}
\label{champs_moderes}

Let us first recall the definition of tame stacks and give a few facts about them. If $\X$ is an Artin stack locally of finite presentation over a scheme $S$, and if its inertia stack $\Ic=\X\times_{\X\times_S\X}\X$ is finite over $\X$, it follows from~\cite{Keel_Mori} that there is a coarse moduli space $\pi : \X \fleche X$ for $\X$. Note by the way that the map $\pi$ is proper. Moreover, the algebraic space $X$ is locally of finite type if $S$ is locally noetherian, and if $\X$ is separated then $X$ is separated as well.

\begin{sousdefi}[\cite{Abramovich_Olsson_Vistoli_Tame_Stacks} 3.1]
\label{def_tame}
 Let $S$ be a scheme and let $\X\fleche S$ be a locally finitely presented algebraic stack over $S$ with finite inertia. Let $\pi : \X \fleche X$ be its Keel-Mori moduli space. We say that $\X$ is tame if the functor $\pi_* : \Qcoh(\X) \fleche \Qcoh(X)$ is exact.
\end{sousdefi}

The tameness condition can be stated in terms of automorphism groups of geometric points. More precisely, a stack $\X \to S$ as above is tame if and only if the automorphism group scheme of any geometric point is linearly reductive (\cite{Abramovich_Olsson_Vistoli_Tame_Stacks}~3.2). If $\X$ is a Deligne-Mumford stack, this means that the order of the automorphism group of any geometric point is prime to the characteristic of the corresponding field. We also recall that the class of tame stacks is stable under arbitrary base change (\cite{Abramovich_Olsson_Vistoli_Tame_Stacks}~3.4) and that if $\X$ is tame then forming its moduli space commutes with any base change (\cite{Abramovich_Olsson_Vistoli_Tame_Stacks}~3.3).

\begin{souslem}
\label{image_directe_plate}
 With the notations and assumptions of~\ref{def_tame}, let $\Fc$ be a quasi-coherent sheaf on $\X$ and let $\Nc$ be a quasi-coherent sheaf on $X$. Then the natural morphism 
$$(\pi_*\Fc)\otimes_{\Oc_{X}}\Nc \flechelongue \pi_*(\Fc\otimes_{\Oc_{\X}}\pi^*\Nc)$$
is an isomorphism. In particular, if $\Fc$ is flat over $S$, then so is $\pi_*\Fc$.
\end{souslem}
\begin{demo}
 We use more or less the same argument as in the proof of~\cite{Abramovich_Olsson_Vistoli_Tame_Stacks}~3.3~(b). Since the question is local on $X$ for the \'etale topology, we can assume that $X$ is an affine scheme. The statement is obvious if $\Nc$ is free. In the general case, let $\Qc_1\fleche \Qc_0 \fleche \Nc \fleche 0$ be a free presentation of $\Nc$. Then we have a commutative diagram with exact rows (since $\pi_*$ is exact):
$$\xymatrix@C=1.9pc{(\pi_*\Fc)\otimes_{\Oc_{X}}\Qc_1 \ar[r] \ar[d] &
(\pi_*\Fc)\otimes_{\Oc_{X}}\Qc_0 \ar[r] \ar[d] & (\pi_*\Fc)\otimes_{\Oc_{X}}\Nc
\ar[r] \ar[d] & 0 \\
\pi_*(\Fc\otimes_{\Oc_{\X}}\pi^*\Qc_1) \ar[r] & \pi_*(\Fc\otimes_{\Oc_{\X}}\pi^*\Qc_0) \ar[r] &
\pi_*(\Fc\otimes_{\Oc_{\X}}\pi^*\Nc) \ar[r] & 0.}$$
The first two columns are isomorphisms, hence so is the third. The last assertion follows immediately, keeping in mind the fact that $\pi_*$ is exact.
\end{demo}

\begin{souslem}[\cite{Mumford_Abelian_Varieties}~\S 5 lemma~1, see also~\cite{EGA} chap.~$0_\textrm{III}$~(11.9.1)]
\label{Mumford_lemma1} $ $
\begin{itemize}
 \item[a)] Let $A$ be a ring and let $C^{\bullet}$ be a complex of $A$-modules such that $C^p\neq 0$ only if $0\leq p\leq n$. Then there exists a complex $K^{\bullet}$ of $A$-modules such that $K^p\neq 0$ only if $0\leq p\leq n$ and $K^p$ is free if $1\leq p\leq n$, and a quasi-isomorphism of complexes $K^{\bullet}\fleche C^{\bullet}$. Moreover, if the $C^p$ are flat, then $K^0$ will be $A$-flat too.
\item[b)] If $A$ is noetherian and if the $H^i(C^{\bullet})$ are finitely generated $A$-modules, then the $K^p$'s can be chosen to be finitely generated.
\end{itemize}
\end{souslem}
\begin{demo}
 The assertion~b) is exactly \cite{Mumford_Abelian_Varieties}~\S 5 lemma~1. For a), the same proof works, erasing the words ``finitely generated'' everywhere.
\end{demo}

\begin{sousprop}
 \label{complexe_coh_tame}
Let $S$ be the spectrum of a ring $A$ (resp. a noetherian ring $A$) and let $\X$ be a quasi-compact and separated (resp. proper) tame stack on $S$. Let $\Fc$ be a quasi-coherent (resp. coherent) sheaf on $\X$ that is flat over $S$. Then there is a complex of flat $A$-modules (resp. of finite type)
$$\xymatrix{0 \ar[r] &M^0 \ar[r] & M^1 \ar[r] & \dots \ar[r] & M^{n} \ar[r] & 0}$$
with $M^i$ free over $A$ for $1\leq i\leq n$, and isomorphisms
$$H^i(M^{\bullet}\otimes_A A') \flechelongue H^i(\X\otimes_A A',\Fc\otimes_A A')\,,\quad i\geq 0$$
functorial in the $A$-algebra $A'$.
\end{sousprop}
\begin{demo}
 Let $\pi : \X \fleche X$ be the moduli space of $X$. Note that $X$ is separated. Choose a finite affine covering $\Ugo=(U_i)_{i\in I}$ of $X$ by affine open subschemes. Then form the \v{C}ech complex $C^{\bullet}=C^{\bullet}(\Ugo,\pi_*\Fc)$ of alternating \v{C}ech cochains. It is a finite complex of flat~(\ref{image_directe_plate}) $A$-modules and it computes the cohomology groups $H^i(X, \pi_*\Fc)$. Since $X$ is separated, the elements of the covering $\Ugo\otimes_A A'$ obtained after a base change $A\fleche A'$ are still affines, so the cohomology of the complex $C^{\bullet}$ is universally isomorphic to the cohomology of $\pi_*\Fc$ on $X$. But this is also the cohomology of $\Fc$ on $\X$, since the functor $\pi_*$ is exact (use for instance the Leray spectral sequence for $\pi$, \cite{Brochard_Picard}~A.2.8). Note that if $A$ is noetherian, $\X$ proper and $\Fc$ coherent, then the modules $H^i(\X, \Fc)$ are finitely generated by \cite{Olsson_lemme_chow}~(1.2), so in this case the cohomology modules of the complex $C^{\bullet}$ are finitely generated. Now it is enough to apply~\ref{Mumford_lemma1} and~\cite{Mumford_Abelian_Varieties}~\S 5 lemma~2.
\end{demo}

\begin{sousremarque}\rm
\label{remarque_coh_tame}
Because of the existence of this complex, all the corollaries that are in~\cite{Mumford_Abelian_Varieties}~\S 5 hold for tame stacks. In other words, the cohomology of such stacks behaves like that of schemes under base change. 
\end{sousremarque}

\subsection{If the base ring has finite global dimension}

First, we prove that in the general case, there is always an \emph{infinite} complex of flat modules computing universally the cohomology of $\Fc$.

\begin{souslem}
 \label{complexe_coh}
 Let $S$ be the spectrum of a ring $A$ and $\X$ a quasi-compact algebraic stack over $S$. Let $\Fc$ be a quasi-coherent sheaf on $\X$. Then there is a complex of $A$-modules
$$\xymatrix{0 \ar[r] &M^0 \ar[r] & M^1 \ar[r] & \dots \ar[r] & M^{n} \ar[r] & \dots }$$
and isomorphisms
$$H^i(M^{\bullet}\otimes_A A') \flechelongue H^i(\X',\Fc')$$
functorial in the $A$-algebra $A'$ (where $\X'=\X\otimes_A A'$ and $\Fc'=\Fc\otimes_A A'$).
If moreover $\Fc$ is flat over $S$, then we can assume that all the $M^i$'s are flat $A$-modules.
\end{souslem}
\begin{demo}
Let $U^0\fleche \X$ be a presentation
of $\X$, such that $U^0$ is an affine scheme. Let $V^1=U^0\times_{\X}U^0$. Let $W^1\fleche V^1$ be a
presentation of the algebraic space $V^1$, the source of which $W^1$ is an affine scheme, and let  $U^1=U^0 \coprod W^1$. We then get a truncated hypercover\footnote{For the definitions of hypercovers, we refer to \cite{SGA4_2} V (7.3.1.2).}
$$\xymatrix{U^1  \ar@<1ex>[r] \ar@<-1ex>[r] & U^0 \ar[l] \ar[r] & \X.}$$
Let $U^{\bullet}$ be the 1-coskeleton\footnote{The 1-coskeleton functor is by definition the right-adjoint of the 1-truncation functor, which to any simplicial object associates its first order truncation.} of this diagram. Clearly this is a
hypercover (of type~1) for $\X$. Moreover,
we can see easily in the construction of the coskeleton (\emph{cf}. \cite{Duskin}~(0.8)) that for every $n\geq 0$, the algebraic stack
$U^{n+2}$ can be expressed in terms of fiber products obtained from the diagram:
$$\xymatrix{U^{n+1} \ar@<1.5ex>[r] \ar@<-1.5ex>[r]^{\vdots}  &U^n.}$$
We deduce that for every $n\geq 0$, $U^n$ is an affine scheme. We denote by $\Fc^i$ the pullback of $\Fc$ on $U^i$. To $U^{\bullet}$ we can associate for every $q$ the alternating chain complex
$$\xymatrix{H^q(U^0,\Fc^0) \ar[r]  & H^q(U^1,\Fc^1) \ar[r] & \dots \ar[r] & H^q(U^p,\Fc^p) \ar[r]&\dots}$$
and we denote by $\check{H}^p(H^q(U^{\bullet},\Fc^{\bullet}))$ the $p$-th cohomology group of this complex. Applying \cite{SGA4_2}~V~(7.4.0.3)
there is a spectral sequence:
$$E_2^{p,q}=\check{H}^p(H^q(U^{\bullet},\Fc^{\bullet}))\Rightarrow H^{p+q}(X,\Fc).$$
Since $\Fc$ is quasi-coherent, we have $H^q(U^i,\Fc^i)=0$ for all $q>0$ and for all $i$, thus $E_2^{p,q}=0$ for all $q>0$. In other words, the spectral sequence degenerates and induces for every $p$ an isomorphism:
$$\xymatrix{\check{H}^p(H^0(U^{\bullet},\Fc^{\bullet})) \ar[r]^-{\sim} & H^p(\X,\Fc).}$$
Now, if $A'$ is an $A$-algebra and $S'=\Spec A'$, the simplicial object $U^{\bullet}\times_S S'$ obtained by base change is an hypercover for $\X'$, and its objects are affine schemes. Thus we also have an isomorphism:
$$\xymatrix{\check{H}^p(H^0(U^{\bullet}\times_S S',\Fc'^{\bullet})) \ar[r]^-{\sim} & H^p(\X',\Fc').}$$
Taking $M^i=H^0(U^i,\Fc^i)$, we get our complex.
Now if $\Fc$ is flat over $S$, the $M^i$'s are obviously flat over $A$.
\end{demo}

The complex given by~\ref{complexe_coh} can be useful in some circumstances. For instance, the following is an immediate corollary.

\begin{souscor}[\cite{Brochard_Picard}~A.3.4]
Let $f : \X \fleche \Y$ be a quasi-compact morphism of $S$-algebraic stacks, and let $\Fc$ be a quasi-coherent sheaf on $\X$. Let $u : \Y' \fleche \Y$ be a flat base change morphism. Let us take the following notations
$$\xymatrix{\X' \ar[r]^v \ar[d]_g \cartesien& \X\ar[d]^f \\
\Y' \ar[r]^u & \Y.}$$
Then for every $q\geq 0$ the natural morphism
$$u^*R^qf_*\Fc \flechelongue (R^qg_*)(v^*\Fc)$$
is an isomorphism.
\end{souscor}

To get deeper results (\emph{e.g.} semicontinuity), we need a finite complex. For that, we would like to truncate the infinite complex given by~\ref{complexe_coh}. So, for a fixed $n$, we want to consider a complex $M'^{\bullet}$ with $M'^i=M^i$ if $i<n$ and $M'^i=0$ for $i>n$. Now we have at least two possibilities for the choice of $M'^n$: either we keep $M'^n=M^n$, but in this case the last cohomology module is changed, or we take $M'^n=\Ker(M^n\fleche M^{n+1})$, but in this case $M'^n$ is not necessarily flat. In both cases, an assumption is missing when we want to replace this finite complex by a finite complex of \emph{finite} modules with the same cohomology (in the first case the last cohomology module is not of finite type and in the second case the last module of the complex might not be flat). In the sequel we will choose the second option, and the whole point is to check that, when the base ring has finite global dimension, we can still replace the complex by a complex of \emph{finite} modules without affecting the first cohomology groups (even after base change!). This is what we do in the following variation of Mumford's lemmas~1 and~2.

\begin{souslem}
Let $n$ be an integer and $A$ a noetherian ring with global cohomological dimension $k \leq n+1$. Let $M^{\bullet}$ be a complex of $A$-modules
$$0\flechelongue M^0 \flechelongue M^1 \flechelongue \dots \flechelongue M^n\flechelongue M^{n+1} \flechelongue 0$$
with $M^i\neq 0$ only if $0\leq i \leq n+1$. Assume that $M^i$ is flat for $0\leq i\leq n$ and that all the cohomology modules $H^i(M^{\bullet})$ are of finite type. Then there exists a complex $K^{\bullet}$ of $A$-modules of finite type and a morphism of complexes $K^{\bullet} \fleche M^{\bullet}$ such that
\begin{itemize}
\item[a)] $K^i\neq 0$ only if $0\leq i\leq n+1$;
\item[b)] $K^i$ is free for $1\leq i \leq n+1$;
\item[c)] $K^0$ is flat;
\item[d)] $K^{\bullet} \fleche M^{\bullet}$ is a quasi-isomorphism;
\item[e)] for any $A$-algebra $B$ and for every $i$, $0\leq i\leq n-k-1$, the morphism 
$$H^i(K^{\bullet}\otimes_A B) \flechelongue H^i(M^{\bullet}\otimes_A B)$$
is an isomorphism.
\end{itemize}
\end{souslem}
\begin{demo}
 Using Mumford's lemma~1~(\cite{Mumford_Abelian_Varieties}~\S 5), there is a complex~$K^{\bullet}$ of finite type $A$-modules and a morphism $f : K^{\bullet} \fleche M^{\bullet}$ satisfying the properties~a), b) and d). The ``mapping cylinder'' $L^{\bullet}$ associated with this morphism of complexes is defined by
\beqn
L^p&=& K^p \oplus M^{p-1} \\
\delta_L(x,y) &=& (\delta_Kx, f(x)-\delta_My) \quad \forall x\in K^p, y\in M^{p-1}
\eeqn
There is a short exact sequence of complexes
$$0\flechelongue M'^{\bullet} \flechelongue L^{\bullet} \flechelongue K^{\bullet} \flechelongue 0$$
where $M'^{\bullet}$ is defined by $M'^p=M^{p-1}$ and $\delta_{M'}=-\delta_M$. This short exact sequence induces a long exact sequence of cohomology
$$\xymatrix@R=0.1pt{0\ar[r] & H^0(L^{\bullet}) \ar[r] & \ar[r]^{\partial} H^0(K^{\bullet}) &\ar[r] H^0(M^{\bullet}) & \\
\ar[r]& H^1(L^{\bullet})\ar[r] &\ar[r]^{\partial} H^1(K^{\bullet}) &\ar[r] H^1(M^{\bullet}) & \\
&\dots &&&}$$
and the cobordism $\partial : H(K^{\bullet}) \fleche H(M^{\bullet})$ coincides with the morphism induced by $f : K^{\bullet} \fleche M^{\bullet}$. Since $f$ is a quasi-isomorphism, this proves that all the $H^i(L^{\bullet})$ vanish, in other words that the complex $L^{\bullet}$ is exact. But $L^{\bullet}$ is:
$$\xymatrix{0 \ar[r]& K^0 \ar[r] & K^1\oplus M^0  \ar[r] & \dots  \ar[r] &K^{n+1}\oplus M^n \ar[r] &M^{n+1} \ar[r] & 0.}$$
All the terms in the middle are flat, thus, splitting this exact sequence in short exact sequences, we get an isomorphism
$$\xymatrix{\Tor_{n+2}^A(M^{n+1},N) \ar[r]^-{\sim}& \Tor_1^A(K^0,N)}$$
for any $A$-module $N$. But the global cohomological dimension of $A$ is $k\leq n+1$, so $\Tor_{n+2}^A(M^{n+1},N)=0$ and $K^0$ is therefore flat.

It remains to prove e). The complex $L^{\bullet}$ is a flat resolution of $M^{n+1}$, so for any $A$-module $N$ and for $i\leq n$, the module $H^i(L^{\bullet}\otimes N)$ is isomorphic to $\Tor_{n+1-i}^A(M^{n+1},N)$. In particular, for $i\leq n-k$ (\emph{i.e.} $n+1-i \geq k+1$) we have~$H^i(L^{\bullet}\otimes N)=0$. Using the ``mapping cylinder'' and the cohomology long exact sequence associated to
$$0\flechelongue M'^{\bullet}\otimes N \flechelongue L^{\bullet}\otimes N \flechelongue K^{\bullet}\otimes N \flechelongue 0$$
we deduce that for any $A$-module $N$ and for every $0\leq i \leq n-k-1$, the natural morphism
$$H^i(K^{\bullet}\otimes_A N) \flechelongue H^i(M^{\bullet}\otimes_A N)$$
is an isomorphism.
\end{demo}

\begin{souscor}
\label{complexe_pour_base_reguliere}
 Let $A$ be a noetherian ring with finite global cohomological dimension, $S$ its spectrum, $\X$ a proper algebraic stack over $S$, and $\Fc$ a coherent $\Oc_{\X}$-module, flat over $S$. Let $n$ be a natural integer. Then there is a finite complex of flat and finite type $A$-modules
$$\xymatrix{0\ar[r]& M^0\ar[r]& M^1\ar[r]& \dots\ar[r]& M^n\ar[r]& M^{n+1},}$$
and functorial isomorphisms
$$\xymatrix{H^i(M^{\bullet}\otimes_A B) \ar[r]^-{\sim}& H^i(\X\otimes_A B, \Fc\otimes_AB)\,,\quad 0\leq i\leq n\,.}$$
\end{souscor}
\begin{demo}
 Let $M^{\bullet}$ be the infinite complex given by the lemma~\ref{complexe_coh}. Since $\X$ is proper and $\Fc$ is coherent, the $A$-modules $H^i(M^{\bullet})$ are of finite type. Let $k$ be the global dimension of $A$ and let $N=n+k+1$. Let $M'^{\bullet}$ be the complex given by
\beqn
M'^p&=& M^p \textrm{ if } p\leq N,\\
M'^{N+1} &=& \Ker(M^{N+1} \fleche M^{N+2}),\\
M'^p &=& 0 \textrm{ if } p\geq N+2.
\eeqn
Apply the previous lemma to $M'^{\bullet}$ and let $K^{\bullet}$ be the resulting complex. Then the $(n+1)$-th truncation of $K^{\bullet}$ is suitable.
\end{demo}

\begin{sousremarque}\rm
 \label{cor_pour_base_reguliere}
As in the case of tame stacks, the existence of this complex implies that all the corollaries from~\cite{Mumford_Abelian_Varieties}~\S 5 hold for a stack over a noetherian ring with finite global dimension (\emph{e.g.} a regular ring with finite dimension).
\end{sousremarque}

\subsection{The semicontinuity theorem}

In this section we will prove the semicontinuity theorem for an arbitrary base ring and for non necessarily tame algebraic stacks. First, let us recall the following fact:

\begin{souslem}
 Let $S$ be a noetherian scheme and $f : S \fleche \N$ a function on $S$. Then $f$ is upper semicontinuous if and only if the two following conditions are satisfied:
\begin{itemize}
 \item[a)] For any discrete valuation ring $A$ and for any morphism $g : \Spec A \fleche S$, we have
$$f(g(\eta))\leq f(g(\xi))$$
where $\eta$ (resp. $\xi$) denotes the generic (resp. special) point of $\Spec A$.
\item[b)] For any noetherian domain $A$ and for any morphism $g : \Spec A \fleche S$, there is a nonempty open subset of $\Spec A$ on which the function $f\circ g$ is constant.
\end{itemize}
\end{souslem}
\begin{demo}
 This is an easy consequence of EGA~$0_{\textrm{III}}$ (\cite{EGA})~9.3.3 and~9.3.4.
\end{demo}

\begin{sousthm}
\label{semicont}
 Let $S$ be a scheme, $\X$ a proper algebraic stack of finite presentation over $S$, and $\Fc$ a coherent $\Oc_{\X}$-module that is flat over $S$. Then for any integer $i\geq 0$, the function
$$\fonction{d_i}{S}{\N}{s}{\dim_{\kappa(s)} H^i(\X_s,\Fc_s)}$$
is upper semicontinuous over $S$.
\end{sousthm}
\begin{demo}
 Obviously we can assume that $S$ is affine, say $S=\Spec A$. By standard limit arguments, we can also assume that $A$ is of finite type over $\Z$. Owing to the previous lemma, it is enough to prove that the theorem holds
\begin{itemize}
 \item[a)] when $A$ is a discrete valuation ring;
\item[b)] over a nonempty open subset of $\Spec A$, when $A$ is a domain.
\end{itemize}
But if $A$ is an integral, finite type $\Z$-algebra, there is a nonempty open subset of $\Spec A$ which is regular. Thus in both cases it is enough to prove the theorem when $A$ is a regular, integral $\Z$-algebra of finite type. Such a ring has finite global cohomological dimension (see \emph{e.g.}~\cite{EGA} chap.~$0_{\textrm{IV}}$~17.3.1). Hence, using the lemma~\ref{complexe_pour_base_reguliere}, there is a finite complex of flat and finite type $A$-modules computing universally the cohomology modules of $\Fc$ over $\X$ at least up to the $i$-th rank. Now, to reach the conclusion that~$d_i$ is upper semicontinuous, we can proceed exactly as in~\cite{Mumford_Abelian_Varieties}~\S 5 (see the corollary p.~50).
\end{demo}

\subsection{Other consequences}
\label{other_csq}

Let us give below some other consequences of the existence of the complexes~\ref{complexe_coh_tame} and~\ref{complexe_pour_base_reguliere}. The results below are generalizations to a stacky context of some standard results in scheme theory. They were needed in the proof of~\ref{thm_ppal}.

\begin{sousprop}
\label{kunneth_tame}
 Let $A$ be a ring with finite global dimension $k$ and let $S$ be its spectrum. Let $f : \X \fleche S$ be a quasi-compact algebraic stack over $S$ and let $\Fc$ be a quasi-coherent sheaf on $\X$ that is flat over $S$. Let $n\in \N$ be an integer. Assume that all the sheaves $R^{n+i}f_*\Fc\, ,\ 1\leq i \leq k$ are flat over $S$. Then forming $R^nf_*\Fc$ commutes with any base change, \emph{i.e.} if $\varphi : S'\fleche S$ is a base change morphism and if $\varphi' : \X\times_S S' \fleche \X$ and $f' : \X\times_S S' \fleche S'$ denote the induced morphisms, then the natural morphism
$$\varphi^*R^nf_*\Fc \flechelongue R^nf'_*(\varphi'^*\Fc)$$
is an isomorphism.
\end{sousprop}
\begin{demo}
 Since forming the higher direct images commutes with any flat base change (\cite{Brochard_Picard}~A.3.4), the assertion is local on $S'$ so we can assume that it is affine, say $S'=\Spec A'$. Since $R^nf_*\Fc$ is the quasi-coherent sheaf corresponding to the $A$-module $H^n(\Xc, \Fc)$ (and similarly over $A'$), it is enough to prove that forming $H^n(\Xc, \Fc)$ commutes with the base change $A\fleche A'$. Let $M^{\bullet}$ be the $(n+k+1)^{\text{th}}$-truncation  (that is, $M^i=0$ if $i\geq n+k+2$) of the complex given by~\ref{complexe_coh}. It is enough to prove that forming $H^n(M^{\bullet})$ commutes with the given base change. But for any $A$-module $N$ there is a K\"unneth spectral sequence
$$E^2_{p,q} = \Tor_{-p}(H^q(M^{\bullet}), N) \Longrightarrow H^{p+q}(M^{\bullet}\otimes_A N).$$
Since all the $H^q(M^{\bullet})$ are flat for $n+1\leq q\leq n+k$, and since the global dimension of $A$ is $k$, the only possibly nonzero terms in this spectral sequence are located at the places marked with a bullet ($\bullet$) in the following picture.

$$\def\latticebody{%
\ifnum\latticeA=1
\else
\ifnum\latticeA=2
\else
  \ifnum\latticeB=2
  \else  
  \ifnum\latticeB=11
  \else  
  \ifnum\latticeB=-1 %
  \else
    \drop{\circ}
  \fi
  \fi
  \fi
  \ifnum\latticeA=-7 \else
  \ifnum\latticeA=-6 \else
    \ifnum\latticeB=0 \drop{\bullet} \fi
    \ifnum\latticeB=1 \drop{\bullet} \fi
    \ifnum\latticeB=3 \drop{\bullet} \fi
    \ifnum\latticeB=9 \drop{\bullet} \fi
  \fi
  \fi
\fi
\fi
\ifnum\latticeA=0
  \ifnum\latticeB=2
  \else  
  \ifnum\latticeB=10
  \else  
  \ifnum\latticeB=-1 %
  \else
    \drop{\bullet}
  \fi
  \fi
  \fi
\fi
}
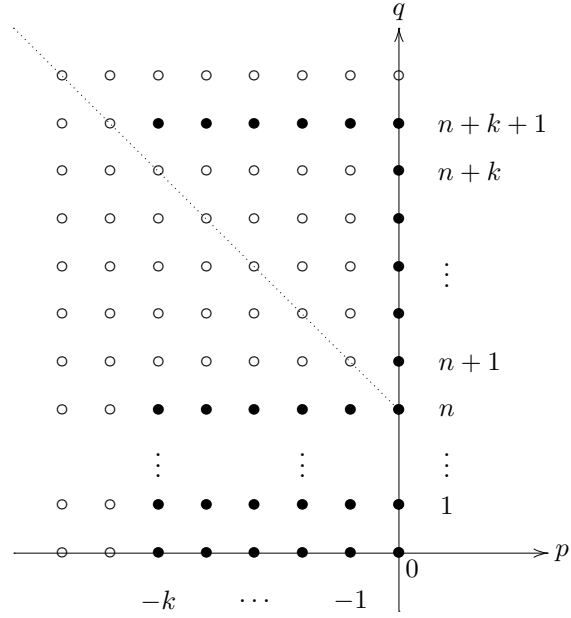
\xy
*\xybox{0;<1.5pc,0mm>:<0mm,1.5pc>::
,0,{\xylattice{-7}2{-1}{10}}
,(0,0)*+<5pt>!U!L{0}
,(1,1)*{1}
,(1,2)*{\vdots}
,(1,3)*{n}
,(1,4)*{\hskip6mm n+1}
,(1,6)*{\vdots}
,(1,8)*{\hskip6mm n+k}
,(1,9)*{\hskip1.2cm n+k+1}
,(-1,-1)*{-1}
,(-3,-1)*{\dots}
,(-5,-1)*{-k}
,(-2,2)*{\vdots}
,(0,3)+LU \ar@{.}
,(-5,2)*{\vdots}
,(-8,11)
}="a"
,{"a"+L \ar "a"+R*+!L{p}}
,{"a"+D \ar "a"+U*+!D{q}}
\endxy
$$
Hence we get an isomorphism $H^n(M^{\bullet})\otimes_A N \fleche H^n(M^{\bullet}\otimes_A N).$
\end{demo}
\begin{sousremarque}\rm
There is also a ``tame version'' of this result, which is left to the reader. 
\end{sousremarque}

\begin{sousprop}
\label{lem_homisom}
Let $\X$ be a proper and flat algebraic stack over a noetherian base scheme $S$ and let $\Lc$, $\Mc$ be invertible sheaves on $\X$. Assume moreover that either $S$ has finite global dimension (\emph{e.g.} $S$ is finite-dimensional and regular) or that $\Xc$ is tame. Then the sheaves $\fHom(\Mc,\Lc) : T \mapsto \Hom_{\Oc_{\X_T}}(\Mc_T,\Lc_T)$ and 
$\fIsom(\Mc,\Lc) : T \mapsto \Isom_{\Oc_{\X_T}}(\Mc_T,\Lc_T)$ are affine schemes of finite type over $S$.
\end{sousprop}
\begin{demo}
 We can assume that $S$ is affine ($S=\Spec A$). Tensoring by $\Mc^{-1}$, we can assume that $\Mc=\Oc_{\X}$. Now the first of the above functors is $T\mapsto H^0(\X_T, \Lc_T)$. Since we already know it is a sheaf, it is enough to consider its restriction to affine schemes, \emph{i.e.} to look at the functor $B \mapsto H^0(\X\otimes_A B, \L\otimes_A B)$. Let $M^{\bullet}$ be a finite complex of finite free $A$-modules computing universally the (beginning of the) cohomology of~$\Lc$. Such a complex is given by~\ref{complexe_coh_tame} or by~\ref{complexe_pour_base_reguliere} depending on our assumption. Then for any $B$, $H^0(\X\otimes_A B, \L\otimes_A B)$ is canonically isomorphic to $H^0(M^{\bullet}\otimes_A B)$. Now, applying~\cite{EGA} EGA~$\textrm{III}_2$~7.4.6 there exists an $A$-module $Q$ of finite type and a functorial isomorphism $H^0(M^{\bullet}\otimes_A B)\fleche \Hom_{A-\textrm{Mod}}(Q,B)$. Thus $\fHom(\Mc,\Lc)$ is the vector bundle $\mathbb{V}(Q)$, hence an affine scheme of finite type. Now $\fIsom(\Mc,\Lc)$ is the closed subscheme of
$\fHom(\Mc,\Lc)\times_S \fHom(\Lc,\Mc)$ defined by the conditions $\varphi\circ \psi=\id$ and $\psi\circ \varphi=\id$ for $(\varphi,\psi) \in \Hom(\Mc_T,\Lc_T)\times \Hom(\Lc_T,\Mc_T)$.
\end{demo}

\begin{sousprop}
\label{images_directes_superieures_affines}
 Let $S$ be a noetherian scheme and let $f : \X \fleche S$ be a proper algebraic stack over $S$. Assume that $S$ has finite global dimension $k$. Let $\Fc$ be a coherent sheaf on $\X$ that is flat over $S$. Assume that all the sheaves $R^{n+i}f_*\Fc$ are flat over $S$ for $0\leq i\leq k$. Then the sheaf $V_n$ on $(\text{\rm Sch}/S)^{\circ}$ defined on affine schemes by
$$V_n(T)=H^n(\X_T, \Fc_T)$$
is an affine scheme of finite type over $S$.
\end{sousprop}
\begin{demo}
We can assume that $S=\Spec A$. Let $M^{\bullet}$ be the complex given by~\ref{complexe_pour_base_reguliere}. Then $V_n$ is (isomorphic to) the Zariski sheaf on $(\textrm{Sch}/S)^{\circ}$ defined by $V_n(A')=H^n(M^{\bullet}\otimes_A A')$ for any $A$-algebra $A'$. Proceeding as in the proof of~\ref{kunneth_tame}, we see that the functor $H^n(M^{\bullet}\otimes_A N)$ in the $A$-module $N$ is isomorphic to $H^n(M^{\bullet})\otimes_A N$, hence it is exact since $H^n(M^{\bullet})$ is flat. Now owing to~\cite{EGA} EGA~$\textrm{III}_2$~7.4.6 there is an $A$-module $Q^n$ and an isomorphism of functors $H^n(M^{\bullet}\otimes_A A') \fleche \Hom_{A-\textrm{Mod}}(Q^n,N)$, thus $V_n$ is the vector bundle $\mathbb{V}(Q^n)$.
\end{demo}

\begin{sousremarque}\rm Once again there is a tame version as follows. Let $\Xc$ be a proper and tame algebraic stack over a noetherian scheme $S$. Let $\Fc$ be a coherent sheaf that is flat over $S$, and such that all the sheaves $R^if_*\Fc$ are flat over $S$ ($i\geq n$). Then the conclusion of the above proposition holds. 
\end{sousremarque}

\paragraph*{Acknowledgments.} I thank Angelo Vistoli, who gave me the idea of using rings with finite global cohomological dimension. I also thank the referee for helpful comments and suggestions.

\bibliographystyle{plain}
\addcontentsline{toc}{section}{References}
\bibliography{mabiblio}

\begin{thebibliography}{10}

\bibitem{SGA6}
{\em Th\'eorie des intersections et th\'eor\`eme de {R}iemann-{R}och}.
\newblock Springer-Verlag, Berlin, 1971.
\newblock S\'eminaire de G\'eom\'etrie Alg\'ebrique du Bois-Marie 1966--1967
  (SGA 6), Dirig\'e par P. Berthelot, A. Grothendieck et L. Illusie. Avec la
  collaboration de D. Ferrand, J. P. Jouanolou, O. Jussila, S. Kleiman, M.
  Raynaud et J. P. Serre, Lecture Notes in Mathematics, Vol. 225.

\bibitem{SGA4_2}
{\em Th\'eorie des topos et cohomologie \'etale des sch\'emas. {T}ome 2}.
\newblock Springer-Verlag, Berlin, 1972.
\newblock S\'eminaire de G\'eom\'etrie Alg\'ebrique du Bois-Marie 1963--1964
  (SGA 4), Dirig\'e par M. Artin, A. Grothendieck et J. L. Verdier. Avec la
  collaboration de N. Bourbaki, P. Deligne et B. Saint-Donat, Lecture Notes in
  Mathematics, Vol. 270.

\bibitem{Abramovich_Olsson_Vistoli_Tame_Stacks}
Dan Abramovich, Martin Olsson, and Angelo Vistoli.
\newblock Tame stacks in positive characteristic.
\newblock {\em Ann. Inst. Fourier (Grenoble)}, 58(4):1057--1091, 2008.

\bibitem{Global_Analysis_1}
Michael Artin.
\newblock Algebraization of formal moduli. {I}.
\newblock In {\em Global Analysis (Papers in Honor of K. Kodaira)}, pages
  21--71. Univ. Tokyo Press, Tokyo, 1969.

\bibitem{BRL}
Siegfried Bosch, Werner L{\"u}tkebohmert, and Michel Raynaud.
\newblock {\em N\'eron models}, volume~21 of {\em Ergebnisse der Mathematik und
  ihrer Grenzgebiete (3) [Results in Mathematics and Related Areas (3)]}.
\newblock Springer-Verlag, Berlin, 1990.

\bibitem{Brochard_qcqsep}
Sylvain Brochard.
\newblock Propri\'et\'es de finitude pour les morphismes non repr\'esentables.
\newblock Unpublished note, available on the webpage of the author.

\bibitem{Brochard_Picard}
Sylvain Brochard.
\newblock Foncteur de {P}icard d'un champ alg\'ebrique.
\newblock {\em Math. Ann.}, 343:541--602, 2009.

\bibitem{Cadman_USTITCOC}
Charles Cadman.
\newblock Using stacks to impose tangency conditions on curves.
\newblock {\em Amer. J. Math.}, 129(2):405--427, 2007.

\bibitem{Duskin}
John~W. Duskin.
\newblock Simplicial methods and the interpretation of ``triple''\ cohomology.
\newblock {\em Mem. Amer. Math. Soc.}, 3(issue 2, 163):v+135, 1975.

\bibitem{EGA}
Alexander Grothendieck.
\newblock \'{E}l\'ements de g\'eom\'etrie alg\'ebrique.
\newblock {\em Inst. Hautes \'Etudes Sci. Publ. Math.}, (4, 8, 11, 17, 20, 24,
  28, 32), 1960-1967.

\bibitem{Kahn_picfini}
Bruno Kahn.
\newblock Sur le groupe des classes d'un sch\'ema arithm\'etique.
\newblock {\em Bull. Soc. Math. France}, 134(3):395--415, 2006.
\newblock With an appendix by Marc Hindry.

\bibitem{Keel_Mori}
Se{\'a}n Keel and Shigefumi Mori.
\newblock Quotients by groupoids.
\newblock {\em Ann. of Math. (2)}, 145(1):193--213, 1997.

\bibitem{LMB}
G{\'e}rard Laumon and Laurent Moret-Bailly.
\newblock {\em Champs alg\'ebriques}, volume~39 of {\em Ergebnisse der
  Mathematik und ihrer Grenzgebiete. 3. Folge. A Series of Modern Surveys in
  Mathematics [Results in Mathematics and Related Areas. 3rd Series. A Series
  of Modern Surveys in Mathematics]}.
\newblock Springer-Verlag, Berlin, 2000.

\bibitem{Mumford_Abelian_Varieties}
David Mumford.
\newblock {\em Abelian varieties}.
\newblock Tata Institute of Fundamental Research Studies in Mathematics, No. 5.
  Published for the Tata Institute of Fundamental Research, Bombay, 1970.

\bibitem{Olsson_lemme_chow}
Martin~C. Olsson.
\newblock On proper coverings of {A}rtin stacks.
\newblock {\em Adv. Math.}, 198(1):93--106, 2005.

\bibitem{Olsson_Hom_stacks}
Martin~C. Olsson.
\newblock {$\underline {\rm Hom}$}-stacks and restriction of scalars.
\newblock {\em Duke Math. J.}, 134(1):139--164, 2006.

\end{thebibliography}
\end{document}